\documentclass[hidelinks,review]{elsarticle}
\usepackage{lineno}
\usepackage{amsmath}
\usepackage{amsfonts}
\usepackage{amsthm}
\usepackage{amssymb}
\usepackage{tikz}
\usepackage{multicol}
\usepackage{multirow}
\usetikzlibrary{fit,positioning,shapes.symbols}

\def\d{{\rm d}}
\def\tE{\text{E}}
\def\T{\text{T}}
\usepackage{subcaption}
\usepackage[top=1.1in, bottom=1.1in, left=1.3in, right=1.3in]{geometry}
\numberwithin{equation}{section}
\numberwithin{figure}{section}
\usepackage{algorithm}
\usepackage{algorithmicx, algpseudocode}
\numberwithin{algorithm}{section}

\newtheorem{theorem}{Theorem}[section]
\newtheorem{lemma}{Lemma}[section]
\newtheorem{remark}{Remark}[section]
\newtheorem{proposition}{Proposition}[section]

\usepackage{mathtools}

\usepackage{tikz}
\usepackage{pgfplots}
\usepgfplotslibrary{polar}
\usepgflibrary{shapes.geometric}
\usetikzlibrary{calc}
\usetikzlibrary{arrows,shapes,intersections}
\usepackage{comment}

 \usepackage{cleveref}
 \crefname{equation}{}{}

\begin{document}
\begin{frontmatter}
\title{On the Application of Stable Generalized Finite Element Method for Quasilinear Elliptic Two-Point BVP}

\author[address1]{T. Aryeni}
\author[address2]{Q. Deng}
\author[address1]{V. Ginting}

\address[address1]{Department of Mathematics \& Statistics, University of Wyoming, Laramie, WY 82071, USA}
\address[address2]{Department of Mathematics, University of Wisconsin-Madison, WI  53706, USA}

\begin{abstract}
In this paper, we discuss the application of the Generalized Finite Element Method (GFEM) to approximate the solutions of quasilinear elliptic equations with multiple interfaces in one dimensional space. The problem is characterized by spatial discontinuity of the elliptic coefficient that depends on the unknown solution. It is known that unless the partition of the domain matches the discontinuity configuration, accuracy of standard finite element techniques significantly deteriorates and standard refinement of the partition may not suffice. The GFEM is a viable alternative to overcome this predicament. It is based on the construction of certain enrichment functions supplied to the standard space that capture effects of the discontinuity. This approach is called stable (SGFEM) if it maintains an optimal rate of convergence and the conditioning of GFEM is not worse than that of the standard FEM. A convergence analysis is derived and performance of the method is illustrated by several numerical examples. Furthermore, it is known that typical global formulations such as FEMs do not enjoy the numerical local conservation property that is crucial in many conservation law-based applications. To remedy this issue, a Lagrange multiplier technique is adopted to enforce the local conservation. A numerical example is given to demonstrate the performance of proposed technique.
\end{abstract}
\end{frontmatter}

\section{Introduction}
\label{intro}
\setcounter{equation}{0}
Realistic mathematical modelings and simulations must often deal with various forms of discontinuities and problems with inherent interfaces. A common illustrative example is simulation of flow and transport of fluids in porous media. An underground soil, for instance, is characterized by spatial heterogeneity, perhaps the simplest one being realized as a layered system, where each layer exhibits a unique conductivity, and with abrupt changes in between. In this context, a relevant motivation comes from modeling the water movement/infiltration in the region near the surface where the pores are filled with water and air (unsaturated zone). Modeling this groundwater flow in unsaturated zone is described by Richards equation 
$$
\partial_t \theta(u) - \nabla \cdot (\kappa({\boldsymbol x},u)\nabla(u - z)) = 0,
$$
which has been proposed by L.A. Richards in 1930 \cite{richards}. The nonlinearity of this equation arises from the water content $\theta $ and the hydraulic conductivity of the soil denoted by $\kappa$ in which both functions are depending on pressure head $u$. The variable $z$ stands for the height against the gravitational direction. The hydraulic conductivity $\kappa$ of different soil types such as sand, clay, and silt has different ability to transmit the water through pore spaces.  This condition will lead to a discontinuity of $\kappa$ at the interfaces of two different types of soil. The situation is made worse by the dependence of $\kappa$ on the pressure head, $u$. Due to the limitation of the availability of the closed-form solution, numerical approximations such as finite difference, finite volume, and finite element are arguably the only reliable procedures to solve this problem (see for example  \cite{doi:10.2136/sssaj2017.02.0058} for a recent review).

Unfortunately, general application of standard finite element approximation to problems of this type fails to maintain an accepted accuracy and convergence optimality. This can be alleviated, for example, by designing discretization of the computational domain to match the discontinuities, which mainly results in a restrictive mesh configuration that cannot allow an interface to cross internal region of the finite element geometry. However, this may not be suitable in many respects. The main reason stems from the lack of knowledge on the exact location of the interfaces. Even when a reasonably adequate information can be gathered on the location of interfaces, it is often realized into an irregular configuration, which in turn presents various challenges in the numerical discretization of the problem. Furthermore, to quantify uncertainty associated with the location of interfaces and the relevant parameters, a common approach is of Monte Carlo type simulations, in which a large number of realizations/samples of interface configurations is proposed and is used as a data, whose results are gathered in the form of some relevant statistics. Since many computational works must be performed, it becomes impractical to change the discretization every time as a different interface configuration is proposed. Against this backdrop is placed an intention to develop a numerical approximation that is flexible toward handling interfaces inherent in the problem. An ideal feature is one that can capture accurately the effects of discontinuity without the necessity to dynamically rearrange discretization of the computational domain. This is indeed desired especially in the realm of Monte Carlo simulations alluded to earlier.

All the above issues give a strong motivation to approximate the solution by a direct extension of the standard FEM called generalized or extended FEM (GFEM/XFEM) that is developed to  handle problems that involve the material discontinuity, moving interface and, crack propagation.  The first development of such methods was recorded in  \cite{babuska97}, and was referred to as Partition of Unity Method (PUM). Later this method was called GFEM in \cite{Copps01,Zhang03,Ivo04,Melenk}. The idea of GFEM is to reduce the discretization errors of the standard FEM without having to change the finite element meshes. This is done by adding more basis that have a compact support around the elements with interface such that they mimic the local behavior of the unknown solution.  

There are several examples of GFEM for the interface problem: Geometric GFEM, Topological GFEM, M-GFEM, and Stable GFEM. The difference between each of these examples is primarily in the way of defining the local enrichment function. Details can be seen in \cite{babuska16,Ivo04}. In topological FEM, the rate of convergence of the $H^1$ semi-norm error is at $h^{1/2}$, where $h$ is mesh parameter, which is similar to the well-known result for standard FEM with uniform mesh and the interface not located at a node. Geometric-FEM and M-FEM have a rate of convergence $h^p$ where $p$ is the degree of Lagrange polynomial interpolation. This rate is similar to the optimal rate of the standard FEM for smooth solution. 

Despite the clear advantage of maintaining optimal convergence properties, studies also found that these GFEMs lead to an extremely high condition number of the stiffness matrix \cite{patrick05,minnebo05}, which gives challenges in solving the corresponding algebraic equations. Later, this issue was addressed in \cite{babuska11} by doing a simple modification on the local enrichment function, creating the Stable GFEM (SGFEM). This reference demonstrates that with this modification, the optimal convergence can be attained without deteriorating the condition number of the system. A further investigation on the conditioning of SGFEM and the comparison with the standard GFEM is given in \cite{babuska16}, confirming that SGFEM maintains the optimal convergences in $H^1$ semi-norm, and the conditioning number of the associated SGFEM matrix is not worse than the standard FEM matrix.

Although a simple modification on the local enrichment function that is suggested in \cite{babuska11} will guarantee a stable GFEM, it is not always the case for some problems in higher dimensions as shown in the application of 2D and 3D fracture mechanics \cite{Gupta13,Gupta15}. Another application for 2D two phase flow problem in \cite{Sauerland13} also shows loss accuracy of the optimal convergence particularly in the case of straight interface problem. Therefore, different modification on the local enrichment is required for a GFEM to be stable.

A lot of studies have made use of the modified enrichment proposed in \cite{babuska11}  and investigated its performance numerically and theoretically for several interface problems. Recently, SGFEM has been implemented in two-dimensional parabolic (time-dependent) interface problem \cite{Zhu19}. The higher order SGFEM has also been developed for the elliptic eigenvalue and source interface problem in \cite{Deng19}. However, there are still not many literatures studying the application of SGFEM for nonlinear interface problems. As reported in \cite{Sauerland13}, SGFEM is applied to industrially relevant two phase/free-surface flow problems governed by the Navier-Stokes equation, which is nonlinear in convection term of its equation. The convergence analysis of the linear SGFEM was established in \cite{babuska11,babuska16,babuska17} and later was generalized to arbitrary order in 1D setting in \cite{Deng19}.

Other relevant motivation of the present investigation is a desire to produce approximate solution that satisfies the local conservation property of the quantity of interest in the presence of the interface system. For Richards' equation in particular, conservative property of the pressure head is needed not only to improve the performance of numerical solutions \cite{Michael90,Kees08}, but also when the resulting velocity $-\kappa({\boldsymbol x},u) \nabla(u-z)$ is coupled to other governing equations, such as for example, concentration of a certain fluid phase invading the soil. In this setting, it is imperative for the approximate velocity to be locally conservative. Several methods such as finite volume method, mixed finite element method, and discontinuous Galerkin method are specifically designed to satisfy this property.  However traditional continuous Galerkin finite element methods fail to yield locally conservative velocity approximation. Applying post-processing technique for such methods has been developed to address this issue. Several work on this subject can be seen in \cite{Deng16,Bernardo07,Bush13,Sun06,Johnson16}. Another technique is called enriched Galerkin (EG) that is done by enriching the approximation space of the CG method with elementwise constant functions \cite{Sun09}. 
Yet another interesting approach was proposed in \cite{Abreu17,Presho15} that is proceeded by constructing the approximate solution that combines the continuous Galerkin formulation and concurrently satisfies the local conservation restrictions. Procedures of this type utilizes a Lagrange multiplier technique, where the approximation is viewed a minimization of the energy functional  over the finite element space under the constraint of algebraic representation of the local conservation property.

In this paper, we investigate an application of high order SGFEM to construct approximate solution of a quasilinear elliptic two-point boundary value problem that possesses a set of discontinuities in its nonlinear elliptic coefficient. This effort can be considered as a first attempt toward the ultimate goal of applying SGFEM to the unsaturated Richards' equation for heterogeneous and layered soil system.  As in \cite{Deng19}, the approximate solution is represented in terms of the usual finite element basis and the enrichment functions aimed at capturing the effects of discontinuity. The resulting nonlinear algebraic system is solved by utilizing Newton's method of iteration. An error analysis in $L^2$ and $H^1$ spaces is carried out that confirms the optimality of SGFEM. Next, we employ the Lagrange multiplier technique as described in \cite{Abreu17} to construct the SGFEM solution that satisfies the local conservation property. This is then validated by a numerical example showing that the optimal convergence behavior of the SGFEM is still maintained and at the same time the local conservation property is satisfied.
 
The outline of this paper is as follows. In \Cref{sec:PS}, we describe the benchmark problem and review the standard continuous Galerkin finite element approximation. In \Cref{sec:SGFEM}, we describe the enriched finite element space that is used in the SGFEM approximation. \Cref{sec:EA} is devoted to a discussion on the existence of the approximate solutions along with convergence and error analysis.  It is then followed by some representative numerical examples in \Cref{sec:NE}. Next in \Cref{sec:FELC}, we present Lagrange multiplier formulation for the FEM/SGFEM solution that satisfies a local conservation property and give a numerical example in \Cref{sec:anelc}. Finally we close the paper with some concluding remarks in \Cref{sec:CO}.
\section{Problem Statement and Standard Finite Element Method}
\setcounter{equation}{0}
\label{sec:PS}
For $K \subset \mathbb{R}$, integer $k\ge 0$ and real number $1\leq p \leq \infty$, we employ standard notation for the Sobolev spaces $W^{k,p}(K)$, with the norm $\Vert\cdot\Vert_{k,p,K}$ and the
seminorm $|\cdot|_{k,p,K}$ \cite{Brenner08,ciarlet2}. In order to simplify the notation, we denote $W^{k,2}(K)$ by $H^k(K)$ and skip the index $p=2$, i.e., $\Vert u \Vert_{k,2,K}= \Vert u \Vert_{k,K}$. We also skip putting $K$ when it is clear that $K$ is the domain of the original problem (later denoted by $\Omega$), thus
we will use $\Vert u\Vert_{k,p,\Omega}=\Vert u \Vert_{k,p},$
$\Vert u \Vert_{k,2,\Omega}= \Vert u \Vert_k$ and $\Vert u\Vert_{0,2,\Omega}=\Vert u \Vert$. The same convention is used for the seminorms as well. In addition, $H_0^1(\Omega)=\{v \in H^1(\Omega)|v=0 \hspace{2mm}\mbox{\rm on}\hspace{2mm} \partial \Omega\}$. In what follows, the symbol $|\cdot |$ will denote the area of a domain, and  $(\cdot,\cdot)$ denote the $L^2(\Omega)$ inner product.

 Let $\Omega = (0,L)\subset\mathbb{R}$, $\Gamma = \{ \gamma_j: j=1,\cdots,m_\Gamma\} \subset \Omega$ with $0<\gamma_j < \gamma_{j+1}<L$ for $j=1,\cdots,m_\Gamma-1$. Let $\Omega_0 = (0,\gamma_1)$, $\Omega_j = (\gamma_j, \gamma_{j+1})$, $j=1,\cdots,m_\Gamma-1$, and $\Omega_{m_\Gamma} = (\gamma_{m_\Gamma}, L)$, such that $\Omega = \cup_{j=0}^{m_\Gamma} \overline{\Omega}_j$. For every $j = 0,\cdots,m_\Gamma$, let $0 < \kappa_{j,\text{min}} \le \kappa_{j,\text{max}}<\infty$ and $\kappa_j: \Omega_j \times \mathbb{R} \to [\kappa_{j,\text{min}}, \kappa_{j,\text{max}}]$ such that $\kappa_j \in C(\overline{\Omega}_j\times \mathbb{R})$ and it is Lipschitz continuous with respect to the second variable, namely,
 \begin{equation} \label{eq:lipschitz}
|\kappa_j(x,\zeta_1) - \kappa_j(x,\zeta_2) | \leq C_{L,j} |\zeta_1 - \zeta_2 |, ~~ \forall  \zeta_1,\zeta_2 \in \mathbb{R}, \text{ and } x \in \Omega_j.
\end{equation}
 Given $f:\Omega \to \mathbb{R}$, the problem is stated as follows: 
 \begin{equation} \label{eq:probstatment}
 \begin{aligned}
&\text{find }  u: \Omega \to \mathbb{R} \text{ governed by }\\
 &\begin{cases}
  \vspace{0.2cm}
 \displaystyle -\frac{\d}{\d x}\left( \kappa_j(x,u(x)) \frac{\d u}{\d x} \right) = f(x), ~~x \in \Omega_j, \text{ for every } j=0,1,\cdots,m_\Gamma,
 \\
 \vspace{0.2cm}
 -\kappa_{j-1}(x,u(x))\dfrac{\d u}{\d x} \Big|_{x = \gamma_j^-} =  -\kappa_{j}(x,u(x))\dfrac{\d u}{\d x} \Big|_{x = \gamma_j^+}, \text{ for every } j=1,\cdots,m_\Gamma,
 \\
u=0 \text{ on } \partial \Omega.
\end{cases}
 \end{aligned}
 \end{equation}

The variational formulation of \cref{eq:probstatment} is to find $u \in H_0^1(\Omega)$ such that
\begin{align} \label{eq:weakform}
a(u;u,w) = \ell(w), ~~\forall w\in H_0^1(\Omega),
\end{align}
where
$$
a(v;u,w):= \sum_{j=0}^{m_\Gamma} \int_{\Omega_j} \kappa_j(x,v) u'(x)w'(x)\d x, \hspace{0.4cm} \ell(w) := \int_\Omega f(x) w(x) \d x.
$$ 
Given $v \in C(\overline\Omega)$,
\begin{equation} \label{eq:bilinearprop}
a(v;w,w) \geq \kappa_{\text{min}} |w|^2_1 ~\text{ and }~ 
a(v;z,w) \leq \kappa_{\text{max}}  |z|_1 \, |w|_1 \le \kappa_{\text{max}}  \|z\|_1 \, \|w\|_1,
\end{equation}
for every $z,w\in H^1(\Omega)$,
where $\displaystyle \kappa_{\text{min}}= \min_{0\le j \le m_\Gamma} \kappa_{j,\text{min}}$ and $\displaystyle \kappa_{\text{max}} = \max_{0\le j \le m_\Gamma} \kappa_{j, \text{max}}$. By Friedrich's inequality, there is a constant $C_0>0$ such that
\begin{equation} \label{eq:coercivity}
C_0 \| w \|_1^2 \le a(v;w,w), ~~\forall w \in H^1_0(\Omega).
\end{equation}
In the forthcoming presentation, finite element approximations are investigated. Existence of such solutions in the appropriate Sobolev space is established, which is then followed by a study of a sequence of the approximations. In particular, it is shown that limit of the sequence satisfies \cref{eq:weakform}, thereby confirming the existence of a weak solution to \cref{eq:probstatment} in $H^1_0(\Omega)$.


%
For the corresponding approximation, we introduce a partition of $\Omega$: 
$0= x_0<x_1<\cdots<x_{N-1}<x_{N} = L$ and set $h_j = x_j - x_{j-1}$ with $h = \max_{1\le j \le N} h_j$.  A restriction that $x_i  \notin \Gamma$ for every $i=0,\cdots,N$ is enforced, which makes a nonconforming partition with respect to $\Gamma$. Denote $\mathcal{T}_h = \{  (x_{j-1},x_j) : j=1,\cdots, N \}$. 
Here it is assumed that an element $\tau \in \mathcal{T}_h$ can contain only one $\gamma \in \Gamma$ or not at all.

The standard continuous finite element space of order $p$ is denoted by
$V_h^p \subset H_0^1(\Omega)$, which contains all continuous piecewise polynomials of degree $p$ vanishing on $\partial \Omega$. Setting $\mathcal{N}_h^p = \{1,2, \cdots, pN-1\}$, this space is characterized as  $V_h^p = \text{span}\{\varphi_j : j\in \mathcal{N}_h^p\}$, where  $\varphi_j$ is the usual nodal Lagrangian polynomial of degree $p$.
The continuous Galerkin finite element approximation to \cref{eq:weakform} reads: find $u_h \in V_h^p$ such that
\begin{equation} \label{eq:weakformfinite}
a(u_h;u_h,w_h) = \ell(w_h), ~~ \forall w_h\in V_h^p.
\end{equation}

Standard practice determines the quality of $u_h \in V_h^p$ through an examination of the approximation property of $V_h^p$, which is usually realized through the interpolation operator $\mathcal{I}^p_h : H^1_0(\Omega) \to V_h^p$ defined as
\begin{equation*}
\mathcal{I}_h^p w = \sum_{j \in \mathcal{N}_h^p} w(t_j) \varphi_j,
\end{equation*}
where $t_j \in \overline{\Omega}$ is such that $\varphi_i(t_j) = \delta_{ij}$. At the elemental level, we may set
\begin{equation} \label{eq:locintp}
\mathcal{I}^p_{\tau} w := \mathcal{I}^p_h w \Big |_\tau = \sum_{j \in \mathcal{N}^p_\tau} w(t_j) \varphi_j,
\end{equation}
where $\mathcal{N}^p_\tau \subset \mathcal{N}^p_h$ is the set of degree of freedom indices associated with a $\tau \in \mathcal{T}_h$. It is also known that the approximation quality of $\mathcal{I}^p_h w$ depends on the smoothness of $w$, in particular (see for example Chapter 1 of \cite{Ern_2004}).
\begin{lemma} \label{lem:insfem}
If $w \in H^1_0(\Omega)$, then $\displaystyle \lim_{h \to 0} | w - \mathcal{I}_h^p w |_1 = 0$. Furthermore, if $w \in H^1_0(\Omega) \cap H^{p+1}(\Omega)$, then $| w - \mathcal{I}_h^p w |_1 \le C h^p | w |_{p+1}$.
\end{lemma}

As indicated in the above lemma, when the function to be approximated is only in $H^1_0(\Omega)$, only convergence is guaranteed;  no information about the optimal convergence order can be gathered.
Raising the quality of the approximation of $V_h^p$ so that this aspect may be displayed hinges on the regularity of the function to be approximated.  However, notice that the presence of the interface system $\Gamma$ in \cref{eq:probstatment} prevents its solution to exhibit a full elliptic regularity. Intuitively, provided that $\kappa_j$ is sufficiently smoother than is prescribed earlier and $f \in H^s(\Omega)$ for $s\ge 0$, then it is expected that any solution of \cref{eq:probstatment} would at most belong to
\begin{equation*}
H^{s+2}_\Gamma(\Omega) := \big \{ w \in H^1_0(\Omega) : w |_{\Omega_j} \in H^{s+2}(\Omega_j), ~\forall j=0,\cdots, m_\Gamma \big \}.
\end{equation*}
 Thus, it is not surprising that when standard continuous Galerkin finite element methods are applied to problems of this type, accuracy of the approximation is suboptimal as hinted in the first part of \Cref{lem:insfem}. Specifically, it will occur when $\mathcal{T}_h$ is not conformed to the interface system $\Gamma$. This gives a motivation for the discussion in the next section.

\section{Stable Generalized Finite Element Methods (SGFEM)}
\setcounter{equation}{0}
\label{sec:SGFEM}
In this section, we adopt the same discretization setting for standard continuous Galerkin finite element method laid out in the previous section. As mentioned earlier, it is assumed that an element $\tau \in \mathcal{T}_h$ can contain only one $\gamma \in \Gamma$. A collection of such elements is
$$
\mathcal{T}_{h,\Gamma} = \{ \tau \in \mathcal{T}_h : \gamma \text{ is located in } \tau, ~\forall \gamma \in \Gamma \}.
$$
The idea of generalized finite element methods is to construct a finite element space that is an enlargement of the standard finite element space by a set of auxiliary functions. These functions are associated with $\tau \in \mathcal{T}_{h,\Gamma}$ such that they capture effects of the discontinuity to the solution.

\begin{figure}[H]
\centering
\includegraphics[scale=0.82]{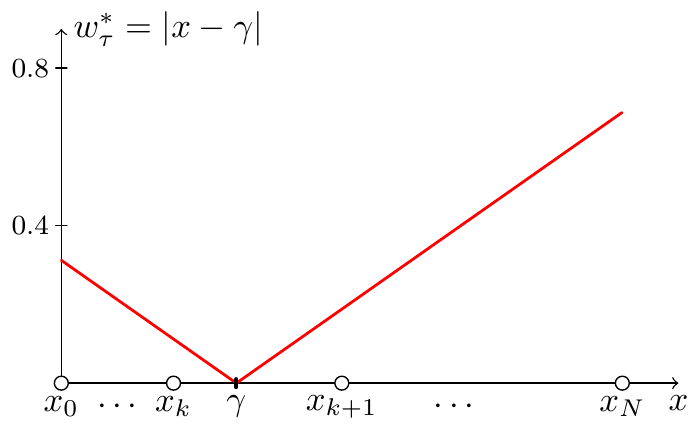}
\hspace*{0.2cm}
\includegraphics[scale=0.82]{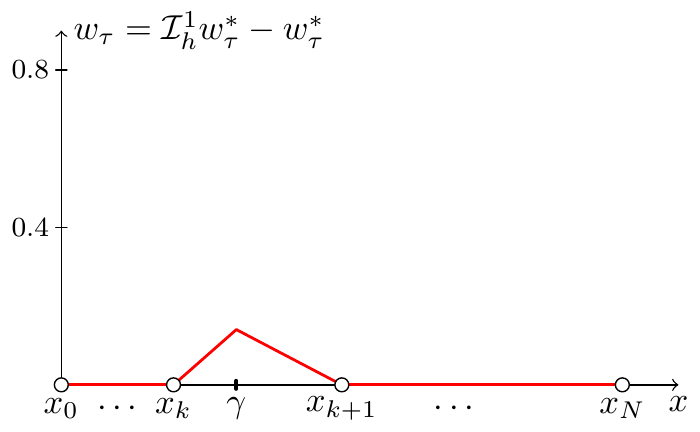}
\caption{An example of $w^*_\tau$ (left plot) and $w_\tau$ (right plot) associated with a $\gamma$ that is located in $\tau = (x_k,x_{k+1}) \in \mathcal{T}_{h,\Gamma} \subset \mathcal{T}_h$}
\label{fig:wwi}
\end{figure}

\begin{figure}[H]
\centering
\includegraphics[scale=0.82]{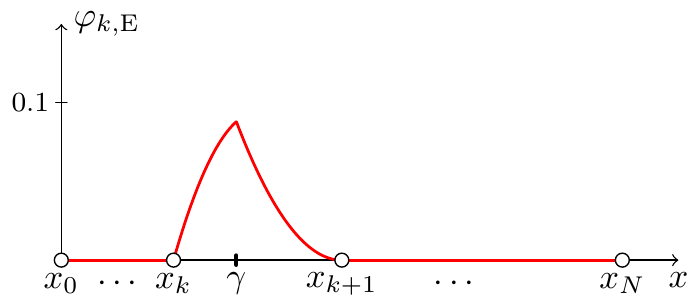}
\hspace*{0.2cm}
\includegraphics[scale=0.82]{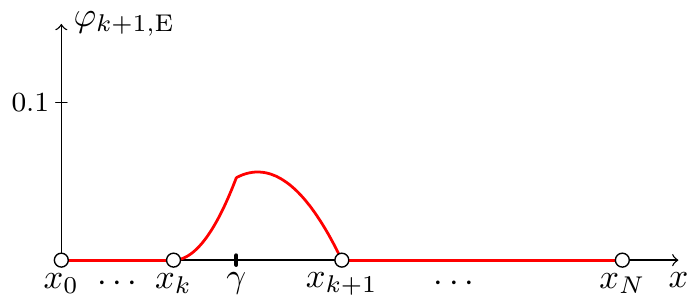}
\caption{The resulting enriched basis functions as applied to $V_h^1$ with $\mathcal{R}_h^1(\tau) = \{ k, k+1 \}$}
\label{fig:ebas}
\end{figure}

Denoting such an enriched space by $V_{h,\tE}^p \subset H_0^1(\Omega)$, it is defined as
$$
V_{h,\tE}^p := V_h^p + V_{h,\text{E}} = \{v_1+ v_2: v_1 \in V_h^p, v_2 \in V_{h,\text{E}}\},
$$
and
$$
V_{h,\text{E}}:= \text{span}\{\varphi_{k,\text{E}}: k \in \mathcal{R}^p_\tau, \forall \tau \in \mathcal{T}_{h,\Gamma}  \},
$$
where $\mathcal{R}^p_\tau \subset \mathcal{N}^p_\tau$ is a set of degree of freedom indices associated with a $\tau \in \mathcal{T}_{h,\Gamma}$. The set $V_{h,\text{E}}$ is called the enrichment space of SGFEM and $\varphi_{k,\text{E}} := w_\tau \varphi_k$, where $\{ w_\tau: \tau \in \mathcal{T}_{h,\Gamma}\}$ is called the set of enrichment functions that are chosen to mimic the true solution near the interfaces. As described in \cite{babuska11,Gupta13,babuska16,babuska17} the choice of $w_\tau$ and $\mathcal{R}^p_\tau$ determines the stability and accuracy properties of the enriched space. For a $\gamma$ that is located in $\tau \in \mathcal{T}_{h,\Gamma}$, the enrichment functions that maintain stability of the approximation are chosen as
\begin{equation} \label{eq:wtau}
w_\tau := \mathcal{I}^1_h w_\tau^* - w^*_\tau, \text{ where } w_\tau^*(x) := |x-\gamma|,
\end{equation}
with $\mathcal{R}^p_\tau = \mathcal{N}^p_\tau$. Obviously $w_\tau$ is piecewise linear and continuous in $\Omega$ with $w_\tau = 0$ outside $\tau$. See \Cref{fig:wwi} for a typical example of $w^*_\tau$ and $w_\tau$ and \Cref{fig:ebas} for the resulting piecewise quadratic enriched basis functions as applied to $V_h^1$.


The stable generalized continuous Galerkin finite element approximation to \cref{eq:weakform} reads: find $u_{h,\tE} \in V_{h,\tE}^p$ such that
\begin{equation} \label{eq:weakformfinitesgfem}
a(u_{h,\tE};u_{h,\tE},w_{h,\tE}) = \ell(w_{h,\tE}), ~~ \forall w_{h,\tE} \in V_{h,\tE}^p.
\end{equation}

Before embarking on the analysis of existence of $u_{h,\tE}$ and its convergence, a discussion on the robustness of $V_{h,\tE}^p$ is warranted, especially on its quality as an approximation space. The following lemma presents the existence of a local interpolation in $V_{h,\tE}^p$ that is applicable to any element with an interface (i.e., any $\tau \in \mathcal{T}_{h,\Gamma}$).  Proof of the lemma below has a slightly different flavor to the one given in \cite{Deng19}.

\begin{lemma} \label{lem:localinterpolation}(local interpolant) Let $\tau = (x_l,x_r) \in \mathcal{T}_{h,\Gamma}$ be associated with a $\gamma \in \Gamma$. Given $v \in C(\overline{\tau})$, let $\mathcal{I}^p_{\tau,{\emph \tE}} v \in C(\overline{\tau})$ such that its restriction to  $[x_l,\gamma]$ and $[\gamma, x_r]$ is a polynomial of degree $p$, and 
$[\mathcal{I}_{\tau, {\emph \tE}} v](\xi_j) = v(\xi_j)$, where $\xi_i = x_l + i(\gamma-x_l)/p$ and $\xi_{p+i+1} = \gamma+i(x_r-\gamma)/p$ for $i = 0,1, \cdots, p$. Then there exists a set of unique $\{ \alpha_j: j=1,\cdots,p+1\} \subset \mathbb{R}$ and $\{ \beta_j: j=1,\cdots,p+1\} \subset \mathbb{R}$ such that
$$
[\mathcal{I}^p_{\tau,{\emph \tE}} v](x) = \sum_{j=1}^{p+1} (\alpha_j + \beta_j w_\tau(x)) \varphi_j(x),
$$
where $\{ \varphi_j: j=1,\cdots,p+1 \}$ is the usual nodal Lagrangian polynomial basis of degree $p$ on $\overline{\tau}$ and $w_\tau$ is as stated in \cref{eq:wtau}.
\end{lemma}
\begin{proof}
Fix a $\gamma$ such that it is located inside $\tau = (x_l,x_r)$.
Given a function $v \in C(\overline{\tau})$, the existence and uniqueness of piecewise polynomial in $\tau$ that interpolates $v$ is clear from the theory of standard nodal polynomial interpolation. In this case, $\mathcal{I}^p_{\tau, \tE} v$ interpolates $v$ at $\{ \xi_i: i=0,\cdots,2p+1\} \subset \mathbb{R}$ as described in the lemma. Set
\begin{equation} \label{eq:interpolating}
Q(x) = \sum_{j=1}^{p+1} (\alpha_j + \beta_j w_\tau(x)) \varphi_j(x), ~~ x\in \overline{\tau},
\end{equation}
where $w_\tau$ is as stated in \cref{eq:wtau}.
The goal is to show the existence of $\{ \alpha_j: j=1,\cdots,p+1\} \subset \mathbb{R}$ and $\{ \beta_j: j=1,\cdots,p+1\} \subset \mathbb{R}$ such that $Q(x) = [\mathcal{I}^p_{\tau,\tE} v](x)$ for every $x \in \overline{\tau}$.

The idea is to construct a linear system governing those coefficients.
Notice that by the above construction, $Q$ is a polynomial of degree at most $p+1$ in $[x_l,\gamma]$ and $[\gamma, x_r]$. Since $\mathcal{I}^p_{\tau,\tE} v$ is a piecewise polynomial of degree $p$, maintaining equality of $Q$ to $\mathcal{I}_{\tau,\tE} v$ requires removing the term $x^{p+1}$ in \cref{eq:interpolating}, yielding an equation
\begin{equation} \label{eq:highestdegree}
\sum_{j=1}^{p+1} a_j \beta_j = 0,
\end{equation}
with $\{ a_j: j=1,\cdots,p+1\} \subset \mathbb{R}$, not all of them are zero. Furthermore, since $\mathcal{I}^p_{\tau,\tE} v$ interpolates $v$ at $2p+1$ distinct points $\xi_i \in \overline{\tau} $, it must satisfy
\begin{equation}\label{eq:interpolationpoints}
Q(\xi_i) = [\mathcal{I}^p_{\tau,\tE} v](\xi_i) = v(\xi_i), ~~i = 0,1, p, p+2,p+3,\cdots, 2p+1,
\end{equation}
where $\xi_{p+1} = \gamma$ has been excluded since interpolatory condition at $\gamma$ has been imposed at $\xi_p = \gamma$. This gives $2p+1$ linear system of equations governing 
$$
\boldsymbol{q} = [\alpha_1, \cdots, \alpha_{p+1},\beta_1,\cdots,\beta_{p+1}]\in \mathbb{R}^{2p+2}.
$$
Combinations of \cref{eq:highestdegree} and \cref{eq:interpolationpoints} gives a $(2p+2)$ linear system
\begin{equation} \label{eq:linsysp}
A \boldsymbol{q} = \boldsymbol{b},
\end{equation}
where  $\boldsymbol{b} = [0, v(\xi_0), v(\xi_1), \cdots, v(\xi_p), v(\xi_{p+2}), v(\xi_{p+3}), \cdots,v(x_{2p+1})] \in \mathbb{R}^{2p+2}$ and $A$ is a square matrix of dimension $2p+2$. If $A$ is nonsingular, then there is a unique $\boldsymbol{q} \in \mathbb{R}^{2p+2}$ satisfying \cref{eq:linsysp}, and the equality of $Q$ to $\mathcal{I}^p_{\tau,\tE} v$ is achieved.

To establish nonsingularity of $A$, it is sufficient to show that $\tilde{\boldsymbol{q}} = \boldsymbol{0} \in \mathbb{R}^{2p+2}$ is the only solution to $A \tilde{\boldsymbol{q}} = \boldsymbol{0}$. But this homogeneous system is equivalent to having
\begin{equation*}
\tilde{Q}(x) = \sum_{j=1}^{p+1} (\tilde{\alpha}_j + \tilde{\beta}_j w_\tau(x) )\varphi_j(x), \text{ with } \tilde{Q}(\xi_i) = 0, \text{ for } i=0,1,\cdots,2p+1, \text{ and }
\sum_{j=1}^{p+1} a_j \tilde{\beta}_j = 0.
\end{equation*}
Thus $\tilde{Q}$ is a piecewise polynomial of degree of at most $p$ having $p+1$ simple zeros in $[x_l,\gamma]$ and $p+1$ simple zeros in $[\gamma,x_r]$. This means it can be expressed as
\begin{equation} \label{eq:roots}
\tilde{Q}(x) = 
\begin{cases}
c_1\displaystyle \prod_{i = 0}^p (x - \xi_i) , ~&x \in [x_l,\gamma],\\
\\
c_2\displaystyle \prod_{i = p+1}^{2p+1} (x-\xi_i), ~&x \in [\gamma,x_r],
\end{cases}
\end{equation}
for some constants $c_1$ and $c_2$.
But, the equation in \cref{eq:roots} gives the leading term $c_1x^{p+1}$ and $c_2x^{p+1}$ that is one degree higher than what is prescribed.  The only way for such a $\tilde{Q}$ to exist is when $c_1 = c_2 = 0$, or equivalently, $\tilde{Q} \equiv 0$. But this implies that every $\tilde{\alpha}_i = 0$ and every $\tilde{\beta}_i = 0$, confirming that $\tilde{\boldsymbol{q}} = \boldsymbol{0}$.

Hence $A$ is invertible, therefore there is a unique $\boldsymbol{q} = [\alpha_1, \cdots, \alpha_{p+1},\beta_1,\cdots,\beta_{p+1}]\in \mathbb{R}^{2p+2}$ satisfying \cref{eq:linsysp} and thus $\mathcal{I}^p_{\tau,\tE} v=Q$ in $\overline{\tau}$. This completes the proof.
\end{proof}
The next lemma is established in \cite{Deng19} for a problem with an interface. We extend it to multiple number of interfaces over $\Omega$.
\begin{lemma} \label{global_int} (global interpolant)
Define $\mathcal{I}^p_{h, {\emph \tE}} : H^1_0(\Omega) \to V^p_{h,{\emph \tE}}$ as follows:
\begin{equation*}
\mathcal{I}^p_{h,{\emph \tE}} v \Big |_\tau =
\begin{cases}
\mathcal{I}^p_{\tau,{\emph \tE}} v &\text{ if } \tau \in \mathcal{T}_{h,\Gamma},\\
\\
\mathcal{I}^p_{\tau} v &\text{ if } \tau \in \mathcal{T}_h \setminus \mathcal{T}_{h,\Gamma},\\
\end{cases}
\end{equation*}
where $\mathcal{I}^p_{\tau,{\emph \tE}}$ is as in \Cref{lem:localinterpolation} and $\mathcal{I}^p_{\tau}$ is as in \cref{eq:locintp}. If $v  \in H^{p+1}_\Gamma(\Omega)$, then
\begin{equation} \label{eq:geint}
\begin{aligned}
|v - \mathcal{I}^p_{h, {\emph \tE}} v|_1 &\leq Ch^p \Bigg(\sum_{j=0}^{m_\Gamma}|v|^2_{p+1,\Omega_j}\Bigg)^{\frac{1}{2}}.
\end{aligned}
\end{equation}
\end{lemma}
\begin{proof}
Fix $v \in H^{p+1}_\Gamma(\Omega)$ (and thus $v \in C(\overline{\Omega}))$. By taking into account the interface system in $\Omega$,
\begin{equation} \label{eq:LLL}
|v - \mathcal{I}^p_{h, \tE} v|^2_1 = \sum_{j=0}^{m_\Gamma} |v - \mathcal{I}^p_{h, \tE} v|^2_{1, \Omega_j},
\end{equation}
so proving \cref{eq:geint} is relegated to establishing a similar estimate for every
$|v - \mathcal{I}^p_{h, \tE} v|_{1, \Omega_j}$.

Since $\mathcal{T}_h$ is nonconforming with respect to $\Gamma$ such that only one $\gamma \in \Gamma$ is located in a $\tau \in \mathcal{T}_{h,\Gamma}$, every $\Omega_j$ has at least one and at most two of such $\tau \in \mathcal{T}_{h,\Gamma}$. Notice that $\mathcal{I}^p_{h,\tE} |_\tau = \mathcal{I}^p_{\tau,\tE}$ for these $\tau$. Furthermore, recall from \Cref{lem:localinterpolation} restriction of $ \mathcal{I}^p_{\tau,\tE}$ to $(x_l,\gamma)$ or $(\gamma,x_r)$ is a polynomial of degree $p$ on $\tau = (x_l,x_r)$. Thus, $\mathcal{I}^p_{h,\tE} v |_{\Omega_j}$ is a piecewise polynomial of degree $p$ that is continuous in $\Omega_j$. Since $v \in H^{p+1}(\Omega_j)$, standard polynomial interpolation estimate (see for example Chapter 1 of \cite{Ern_2004}) gives
$$
 |v - \mathcal{I}^p_{h, \tE} v|_{1, \Omega_j} \le C h^p |v|_{p+1,\Omega_j},
$$
which on its substitution in \cref{eq:LLL} results in \cref{eq:geint}.
\end{proof}

To reiterate what has been mentioned earlier, when the function to be interpolated does not enjoy a global smoothness (in this case in $H^{p+1}(\Omega)$) due to the presence of $\Gamma$ and $\mathcal{T}_h$ does not conform with $\Gamma$, then the standard interpolation operator $\mathcal{I}^p_h$ cannot attain the optimal approximation property of $\mathcal{I}^p_{h, {\tE}}$. This is translated into the performance of the approximate solutions, which will be made clear in the error analysis. In a related matter, the next lemma is particularly needed in the error analysis to handle the quasilinear nature of the original problem, i.e., due to the appearance of $\kappa_j(x,u)$ in \cref{eq:probstatment}.

\begin{lemma}\label{lem:norm16}
If $v  \in H^{p+1}_\Gamma(\Omega)$, then
\begin{equation} \label{eq:geint16}
\begin{aligned}
|v - \mathcal{I}^p_{h, {\emph \tE}} v|_{1,6} &\leq h^{p-\frac{1}{3}} \Bigg(\sum_{j=0}^{m_\Gamma}|v|^2_{p+1,\Omega_j}\Bigg)^{\frac{1}{2}},
\end{aligned}
\end{equation}
where $\mathcal{I}^p_{h, {\emph \tE}}$ is the global interpolant defined in \Cref{global_int}.
\end{lemma}
\begin{proof}
Similar to the previous lemma, we may write
\begin{align} \label{eq:sumj}
|v - \mathcal{I}^p_{h,\tE} v|_{1,6} =\Bigg( \sum_{j=0}^{m_\Gamma} |v -\mathcal{I}_{h,\tE}^p v |^6_{1,6,\Omega_j} \Bigg)^{\frac{1}{6}}.
\end{align}
Because $\mathcal{I}^p_{h,E} v|_{\Omega_j}$ is a continuous piecewise polynomial of degree $p$, we can further write
\begin{align} \label{eq:normsum}
|v - \mathcal{I}_{h,\tE}^p v|^6_{1,6,\Omega_j} = \sum_{\tau \in \mathcal{T}_j} | v - \mathcal{I}_{\tau}^p v|^6_{1,6,\tau} + \sum_{\tau_{1/2} \in \mathcal{T}_{j,\Gamma}}  | v - \mathcal{I}_{\tau,\tE}^p v|^6_{1,6,\tau_{1/2}},
\end{align}
where $\mathcal{T}_j = \{ \tau \in  \mathcal{T}_h \setminus \mathcal{T}_{h,\Gamma} :   \tau \cap  \Omega_j \ne \varnothing  \}$ and  $\mathcal{T}_{j,\Gamma} = \{ \tau \cap \Omega_j : \tau \in \mathcal{T}_{h,\Gamma}  \}$, and thus the union of all these intervals is $\Omega_j$.
The first summation is coming from adding the interpolation error over all elements $\tau$ in $\Omega_j$ that have no interface $\gamma \in \Gamma$, whereas the second summation is the interpolation error in half portion of $\tau \in\mathcal{T}_{h,\Gamma}$ belonging in $\Omega_j$ and is denoted by $\tau_{1/2} = \tau \cap \Omega_j$. Note that for every $\Omega_j$, there are at most two of such $\tau_{1/2}$. Recall that both $\mathcal{I}_{\tau}^p v |_\tau$ for any $\tau \in \mathcal{T}_j$
and $\mathcal{I}_{\tau,\tE}^p v |_\tau$ for any $\tau_{1/2} \in \mathcal{T}_{j,\Gamma}$ are polynomials of degree $p$ that interpolate $v$ in $\tau$.

Let $\varsigma$ be either $\tau \in \mathcal{T}_j$ or $\tau_{1/2} \in \mathcal{T}_{j,\Gamma}$, and define $e : \varsigma \rightarrow \mathbb{R}$ as $e =  \big(v - \mathcal{I}_{h,\tE}^p v\big)|_{\varsigma}$.
Applying the fundamental theorem of calculus to the $k^{th}$  order derivative of $e $, $k \leq p$, we have for any point $x \in  \overline{\varsigma}$,
\begin{align}  \label{eq:ftc}
e^{(k)}(x) = e^{(k)}(\bar{x}_k) + \int_{\bar{x}_k}^x e^{(k+1)}(t_k)\ \d t_k = \int_{\bar{x}_k}^x e^{(k+1)}(t_k) \ \d t_k, 
\end{align} 
where $\bar{x}_k \in \varsigma$ with $e^{(k)}(\bar{x}_k) = 0$ (guaranteed by Rolle's theorem). 
Applying \eqref{eq:ftc} up to $p^{th}$ derivative of $e$, we can express the first derivative of $e$ as follows
\begin{align*}
e^\prime(x)  =  \int_{\bar{x}_1}^x e^{\prime\prime}(t_1)\ \d t_1 = \int_{\bar{x}_1}^x \int_{\bar{x}_2}^{t_1} e^{\prime \prime \prime}(t_2) \ \d t_2 \d t_1=  \int_{\bar{x}_1}^x \int_{\bar{x}_2}^{t_1} \cdots \int_{\bar{x}_p}^{t_{p-1}} e^{(p+1)}(t_p) \ \d t_{p} \d t_{p-1} \cdots \d t_1.
\end{align*}
Note that $(\mathcal{I}_{h,E}^pv)^{(p+1)}(x) = 0$ for every $x \in \overline{\varsigma}$, so $e^{(p+1)} = v^{(p+1)}|_\varsigma$, which can then be used to bound the following integral along with an application of Cauchy-Schwarz inequality:
\begin{align*}
 \int_{\bar{x}_p}^{t_{p-1}} e^{(p+1)}(t_p) \ \d t_{p} \leq \int_{\varsigma} |v^{(p+1)}(t_p)| \ \d t_{p} \leq h^{\frac{1}{2}} \|v^{(p+1)}\|_{0,\varsigma},
\end{align*}
where $h_\varsigma$ has been bounded by $h$. Thus,
\begin{align*}
|e'(x)| &\leq  \int_{\bar{x}_1}^x \int_{\bar{x}_2}^{t_1} \cdots \int_{\bar{x}_{p-1}}^{t_{p-2}} h^{\frac{1}{2}} \|v^{(p+1)}\|_{0,\varsigma}  \ \d t_{p-1} \cdots d_{t_1} \\
&\leq h^{\frac{1}{2}} \|v^{(p+1)}\|_{0,\varsigma} h^{p-1} = h^{p-\frac{1}{2}}  \|v^{(p+1)}\|_{0,\varsigma}.
\end{align*}
Taking the $L^6$-norm of $e^\prime$ over $\varsigma$ gives
\begin{equation*}
|e|^6_{1,6,\varsigma} = \int_\varsigma |e'(x)|^6 \ \d x \le \int_\varsigma h^{6(p-\frac{1}{2})} \|v^{(p+1)}\|^6_{0,\varsigma}\ \d x = h^{6p-2} \|v^{(p+1)}\|^6_{0,\varsigma}.
\end{equation*}
Substituting this into \eqref{eq:normsum} yields
\begin{align*}
|v - \mathcal{I}_{h,\tE}^p v|^6_{1,6,\Omega_j} \leq \sum_{\tau \in \mathcal{T}_j} h^{6p-2} \|v^{(p+1)}\|^6_{0,\tau} +  \sum_{\tau_{1/2} \in \mathcal{T}_{j,\Gamma}} h^{6p-2} \|v^{(p+1)}\|^6_{0,\tau_{1/2}}.
\end{align*}
Because $v \in  H^{p+1}(\Omega_j)$, we can combine the summation, therefore
\begin{align*}
|v - \mathcal{I}_{h,\tE}^p v|^6_{1,6,\Omega_j} \leq h^{6p-2} \|v^{(p+1)}\|^6_{0,\Omega_j}.
\end{align*}
Upon substitution of this last inequality to \eqref{eq:sumj} yields
$$
\begin{aligned}
|v - \mathcal{I}^p_{h,\tE} v|_{1,6} &\le \Bigg( \sum_{j=0}^{m_\Gamma} h^{6p-2} \|v^{(p+1)}\|^6_{0,\Omega_j} \Bigg)^{\frac{1}{6}}\\
&= h^{p-\frac{1}{3}} \Bigg( \Bigg( \sum_{j=0}^{m_\Gamma}  \|v^{(p+1)}\|^6_{0,\Omega_j} \Bigg)^{\frac{1}{3}} \Bigg)^{\frac{1}{2}}\\
&\le h^{p-\frac{1}{3}} \Bigg( \sum_{j=0}^{m_\Gamma}  \|v^{(p+1)}\|^2_{0,\Omega_j} \Bigg)^{\frac{1}{2}},
\end{aligned}
$$
and the proof is complete.
\end{proof}

\begin{remark}
\Cref{lem:norm16} can be generalized for any positive integer $r \in [2,\infty)$, so that 
$$
|v - \mathcal{I}^p_{h, {\emph \tE}} v|_{1,r} \le \mathcal{O} ( h^{p -\frac{1}{2} + \frac{1}{r}} ).
$$
\end{remark}

\section{An Analysis}
\setcounter{equation}{0}
\label{sec:EA}
 In this section, we give an analysis pertaining to the approximations of the solution of \cref{eq:probstatment}. We begin with establishing existence of the approximation and demonstrate that it converges to a weak solution of \cref{eq:probstatment}. This is then followed by an error estimation of the SGFEM solution. Various  mathematical tools and techniques used in the analysis can be seen for example in \cite{MR431747} and \cite{MR1275952}.  Due to procedural similarity in conducting the analysis, in what follows, the finite element space in which the approximation is sought is generically denoted by $\widetilde{V}_h \subset H^1_0(\Omega)$, where $\widetilde{V}_h$ is either $V_h^p$ for standard continuous Galerkin FEM or $V_{h,\tE}^p$ for SGFEM. 

 \subsection{Existence of Approximate Solutions and Convergence Analysis}
Existence of the approximate solutions and their convergence require an assumption that there is $\widetilde{\mathcal{I}}_h : H^1_0(\Omega) \to \widetilde{V}_h$ such that 
\begin{equation} \label{eq:Itilde}
\| \widetilde{\mathcal{I}}_h v - v \|_1 \to 0 \text{ as } h \to 0.
\end{equation}
An example of such an operator is established for $\widetilde{V}_h = V_h^p$ by utilizing 
\Cref{lem:insfem} and for $\widetilde{V}_h = V_{h,\tE}^p$ by utilizing \Cref{global_int} with the help of Friedrich's inequality.

\begin{theorem} \label{thm:41}
If $f \in L^2(\Omega)$, then there exists a $\widetilde{u}_h \in \widetilde{V}_h$ governed by
\begin{equation} \label{eq:finitevariational}
a(\widetilde{u}_h;\widetilde{u}_h,\widetilde{v}_h) = \ell(\widetilde{v}_h), ~~ \forall \widetilde{v}_h \in \widetilde{V}_h.
\end{equation}
\end{theorem}

\begin{proof}
Consider a mapping $\T: \widetilde{V}_h \rightarrow \widetilde{V}_h$ defined by the relation
\begin{equation} \label{eq:Th}
a(y;\T(y),\widetilde{v}_h) = \ell(\widetilde{v}_h), ~~\forall \widetilde{v}_h \in \widetilde{V}_h.
\end{equation}
In this regard, existence of $\widetilde{u}_h \in \widetilde{V}_h$ satisfying \eqref{eq:finitevariational} is equivalent to showing that $\T$ has a fixed point in $\widetilde{V}_h$.
For a given $y \in \widetilde{V}_h$, existence of a unique $\T(y)$ is established by the Lax-Milgram theorem (see for example, p. 317 of \cite{Evans10}).
By setting $\widetilde{v}_h = \T(y)$ in \eqref{eq:Th},  and using \eqref{eq:coercivity} and the boundedness of $\ell$, we get
\begin{align} \label{eq:boundforT}
C_0 \|\T(y)\|_1^2 \leq a(y;\T(y),\T(y) )= \ell(\T(y) )\leq \| f \| \, \|\T(y)\|_1,
\end{align}
from which we confirm that $\T(y) \in \widetilde{V}_h$ is bounded, i.e.,
\begin{equation} \label{eq:mapping}
\|\T(y)\|_1 \le C_0^{-1} \| f \|.
\end{equation}
Next, to show the continuity of $\T$, it is sufficient to demonstrate that it is Lipschitz continuous. Replacing $\widetilde{v}_h$ in \eqref{eq:Th} by $\T(y) - \T(z)$, and using  \eqref{eq:coercivity} and linearity of $a(\cdot ; \cdot, \cdot)$ on the second argument yields
\begin{align*}
C_0 \|\widetilde{v}_h\|_1^2 &\leq  a(y;\T(y),\widetilde{v}_h) - a(y;\T(z),\widetilde{v}_h)\\
&= \ell(\widetilde{v}_h) - a(z;\T(z),\widetilde{v}_h) + a(z;\T(z),\widetilde{v}_h) -  a(y;\T(z),\widetilde{v}_h)\\
&= \ell(\widetilde{v}_h) - \ell(\widetilde{v}_h) +  a(z;\T(z),\widetilde{v}_h) -  a(y;\T(z),\widetilde{v}_h)\\
&= a(z;\T(z),\widetilde{v}_h) -  a(y;\T(z),\widetilde{v}_h).
\end{align*}
By Lipschitz continuity of $\kappa_j$ and Cauchy-Schwarz inequality,
\begin{align*}
a(z;\T(z),\widetilde{v}_h) -  a(y;\T(z),\widetilde{v}_h) &\leq \sum_{j=0}^m \int_{\Omega_j}|\kappa_j(x,z) - \kappa_j(x,y)| \, |[\T(z)]'(x)| \, |\widetilde{v}_h^\prime(x)| \, \d x \\
&\leq \sum_{j=0}^m \int_{\Omega_j} C_{L,j} |z(x) -y(x)| \, |[\T(z)]'(x)| \, |\widetilde{v}_h^\prime(x) | \, \d x  \\
&\leq C_L \int_{\Omega} |z(x) -y(x)| \, |[\T(z)]'(x)| \, |\widetilde{v}_h^\prime(x)| \, \d x \\
&\leq C_L \|y - z\|_{0,\infty}\|\T(z)\|_1 \|\widetilde{v}_h\|_1, 
\end{align*}
where $\displaystyle C_L = \max_{0 \le j \le m} C_{L,j}$. Note that since $y,z \in \widetilde{V}_h \subset H^1_0(\Omega)$,  $\|y-z\|_{0,\infty} \leq C_\Omega |y-z|_1$, which together with \eqref{eq:mapping} implies 
$$
C_0 \|\widetilde{v}_h\|_1 \leq C_L C_0^{-1} C_\Omega  \| f \| \, \|y- z\|_1.
$$
Thus $\T$ is Lipschitz continuous i.e., $\|\T(y) -\T(z)\|_1 \leq C_L C_0^{-2} C_\Omega \| f \| \, \|y - z\|_1$. Since $\T$ is continuous, existence of a $\widetilde{u}_h \in \widetilde{V}_h$ satisfying $\T(\widetilde{u}_h) = \widetilde{u}_h$ is guaranteed by the Brouwer Fixed Point Theorem. This completes the proof.
\end{proof}
\begin{remark}
 Note that the Brouwer Fixed Point Theorem does not guarantee the uniqueness of $\widetilde{u}_h$. If in addition $f \in L^2(\Omega)$ is chosen such that $C_L C_0^{-2} C_\Omega \| f \| < 1$, then $\text{\em T}$ in the above proof is actually a contraction. In this setting, existence and uniqueness of $\widetilde{u}_h$ can be obtained from the Banach Fixed Point Theorem.
\end{remark}

In the following theorem, we show the existence of a weak solution of \eqref{eq:probstatment} as a weak limit of the Galerkin approximations $\widetilde{u}_h \in \widetilde{V}_h$. 

\begin{theorem} \label{thm:existence222}
Let $\{\widetilde{V}_h\}_{h\rightarrow 0}$ be a family of finite dimensional subspaces of $H_0^1(\Omega)$ and let  $\{\widetilde{u}_h\}_{h\rightarrow 0} $ be a sequence of the Galerkin approximations satisfying \eqref{eq:finitevariational}. Then
there exists a subsequence $\{\widetilde{u}_{\bar{h}}\} \subset \{\widetilde{u}_h\}$ and an element $u \in H^1_0(\Omega)$ such that
\begin{equation} \label{eq:weakcon}
\widetilde{u}_{\bar{h}} \rightharpoonup u \in H^1_0(\Omega) \text{ as } \bar{h} \rightarrow 0,
\end{equation}
and $u$ is a weak solution of  \eqref{eq:probstatment}, i.e., it satisfies $a(u; u,w) = \ell(w)$ for every $w \in H^1_0(\Omega)$.
\end{theorem}

\begin{proof}
First, existence of $\{\widetilde{u}_h\}$ satisfying \eqref{eq:finitevariational} is already established in \Cref{thm:41}. Furthermore, by \eqref{eq:mapping}, $\|\widetilde{u}_h\|_1 \leq C_0^{-1} \| f \|$.
Since  $\{\widetilde{u}_h\}$ is bounded in $H^1_0(\Omega)$, it has a subsequence $\{\widetilde{u}_{\bar{h}}\} \subset \{\widetilde{u}_h\}$ that is converging weakly to a limit in $H_0^1(\Omega)$ (see for example, p. 726 of \cite{Evans10}). Suppose $u \in H_0^1(\Omega)$ is the weak limit of subsequence $\{\widetilde{u}_{\bar{h}}\}$ such that \eqref{eq:weakcon} holds, then for any $\mathcal{L} \in [H_0^1(\Omega)]^\ast$,
\begin{equation} \label{eq:weakconvergencecon}
\mathcal{L}(\widetilde{u}_{\bar{h}}) \rightarrow \mathcal{L}(u) \text{   as } \bar{h} \rightarrow 0.
\end{equation}
Furthermore, the  Rellich-Kondrachov theorem in (see for example, p. 288 of \cite{Evans10}) says that the subsequence $\{\widetilde{u}_{\bar{h}}\}$ converges strongly to $u$ in $L^2(\Omega)$, i.e,
\begin{equation} \label{eq:strongcon}
 \|\widetilde{u}_{\bar{h}} - u\| \to 0 \text{   as } {\bar{h} \rightarrow 0}.
\end{equation}
Now we show that $u$ is governed by $a(u;u,w) = \ell(w)$ for every $w\in H_0^1(\Omega)$.
Consider an arbitrary $v \in  C_0^\infty(\Omega)$  and let $\widetilde{\mathcal{I}}_h v  \in \widetilde{V}_h$ be its approximation that satisfies \eqref{eq:Itilde}. Using \cref{eq:finitevariational}, adding and subtracting $a(u;\widetilde{u}_{\bar{h}},v)$ and $a(\widetilde{u}_{\bar{h}};\widetilde{u}_{\bar{h}},v)$,
\begin{equation} \label{eq:auvlv}
|a(u;u,v) - \ell(v)| =  |a(u;u,v) - a(\widetilde{u}_{\bar{h}};\widetilde{u}_{\bar{h}},\widetilde{\mathcal{I}}_h v) +\ell(\widetilde{\mathcal{I}}_h v) - \ell(v)| \le I_1 + I_2 + I_3 + I_4,
\end{equation}
where
$$
\begin{aligned}
I_1 &= |a(u;u,v) - a(u;\widetilde{u}_{\bar{h}},v)|,\\
I_2 &= |a(u;\widetilde{u}_{\bar{h}},v) - a(\widetilde{u}_{\bar{h}};\widetilde{u}_{\bar{h}},v)|,\\
I_3 &= |a(\widetilde{u}_{\bar{h}};\widetilde{u}_{\bar{h}},v) - a(\widetilde{u}_{\bar{h}};\widetilde{u}_{\bar{h}},\widetilde{\mathcal{I}}_h v|,\\
I_4 &= |\ell(\widetilde{\mathcal{I}}_h v) - \ell(v)|.
\end{aligned}
$$

Since  $a(\cdot; \cdot,\cdot)$ is bounded in $H^1_0(\Omega)$, then $\mathcal{L}(w) = a(u; w,v)$ is bounded, and by the fact that $u_{\bar{h}}  \rightharpoonup u$ in $H^1_0(\Omega)$, it is clear that
$a(u;\widetilde{u}_{\bar{h}},v) \rightarrow a(u;u,v)$ as  $\bar{h} \rightarrow 0$, resulting in $I_1 \to 0$.
 
 Taking into account the Lipschitz continuity of  $\kappa_j$ yields the following estimate
\begin{align*}
I_2 &\leq  \sum_{j=0}^m \int_{\Omega_j} |\kappa_j(x,u) - \kappa_j(x,\widetilde{u}_{\bar{h}})| \, |\widetilde{u}_h^\prime(x)| \, |v^\prime(x)| \, \d x \\
&\leq \sum_{j=0}^m C_{L,j} \int_{\Omega_j} |u(x) - \widetilde{u}_{\bar{h}}(x)| \, |\widetilde{u}_{\bar{h}}^\prime(x)| \, |v^\prime(x)| \, \d x \\
& \leq C_L \|v^\prime\|_{0,\infty} \int_\Omega  |u(x) - \widetilde{u}_{\bar{h}}(x)| \, |\widetilde{u}_{\bar{h}}^\prime(x)| \, \d x \\
&\le C_L   \|v^\prime\|_{0,\infty} \, \|u - \widetilde{u}_{\bar{h}}\| \, \|\widetilde{u}_{\bar{h}} \|_1\\
&\le C_L  C_0^{-1} \| f \| \, \|v^\prime\|_{0,\infty} \, \|u - \widetilde{u}_{\bar{h}} \|.
\end{align*}  
 Utilizing  \cref{eq:strongcon}, it is confirmed that $I_2 \to 0$ as $\bar{h} \to 0$.
 
Due to the boundedness of $a(\cdot;\cdot,\cdot)$, $I_3 \leq \kappa_{\text{max}} \| u_h \|_1 \|v - \widetilde{\mathcal{I}}_h v \|_1$, which along with \cref{eq:Itilde} establishes $I_3 \to 0$ as $h \to 0$. Likewise, 
$I_4 \leq \| f \| \|\widetilde{\mathcal{I}}_h v - v\|_1 \rightarrow 0$ as $h \to 0$.

By taking into consideration convergence of all these terms back in \cref{eq:auvlv}, we arrive at
\begin{equation} \label{eq:varform}
a(u;u,v) = \ell(v), \hspace{0.5cm} \forall v\in C_0^\infty(\Omega).
\end{equation}

Finally, recall that $ C_0^\infty(\Omega) $ is dense in $H_0^1(\Omega)$, so that given $w \in H_0^1(\Omega)$ there exists a sequence $(v_i) \subset  C_0^\infty(\Omega)$ such that 
$\|w - v_i \|_1 \to 0$ as $i \to \infty$. By using $v = v_i$ in \eqref{eq:varform} and the boundedness of $a(\cdot ; \cdot, \cdot)$ and $\ell(\cdot)$, 
\begin{align*}
|a(u;u,w) - \ell(w) | &\le |a(u;u,w-v_i) | +  |\ell(v_i - w)|
\leq (\kappa_{\text{max}} \|u\|_1 + \| f \|) \|w - v_i\|_1 \to 0,
\end{align*}
as $i \to \infty$.
Therefore $u \in H_0^1(\Omega)$ is a weak solution of \eqref{eq:probstatment}.
\end{proof}

\subsection{An Error Analysis for the SGFEM Solution}
In this section, a detailed error analysis of the SGFEM solution is presented. The main purpose behind the analysis is to demonstrate that under the assumption that the solution of \cref{eq:probstatment} belongs to
$H^{p+1}_\Gamma(\Omega)$, then its approximation sought in $V_{h,\tE}^p$ maintains the convergence optimality. As stated earlier, this is a desirable trait that the standard finite element space $V_{h}^p$ cannot achieve when $\mathcal{T}_h$ is not conformed to the interface system $\Gamma$. As before, let $u \in H^1_0(\Omega)$ be a weak solution of  \cref{eq:probstatment} and let $\widetilde{u}_h \in \widetilde{V}_h$ be its approximation, which is governed by \cref{eq:finitevariational}.

\begin{lemma}\label{H1ineq}
There exists a positive constant $C$ independent of $h$ and $u$ such that
\begin{equation} \label{eq:sebutsajadia}
|u - \widetilde{u}_h|_1 \leq C( |u - \widetilde{w}_h |_1 + | \widetilde{w}_h |^2_{1,6} \|u - \widetilde{u}_h\|), ~\forall ~\widetilde{w}_h \in \widetilde{V}_h.
\end{equation}
\end{lemma}
\begin{proof}
Given any $\widetilde{w}_h \in \widetilde{V}_h$, triangle inequality gives
\begin{equation} \label{eq:triangle}
|u - \widetilde{u}_h|_1 \leq |u - \widetilde{w}_h|_1 + |\widetilde{w}_h -\widetilde{u}_h|_1,
\end{equation}
so the remainder of the proof is concentrated on estimating $\widetilde{e}_h = (\widetilde{w}_h -\widetilde{u}_h) \in \widetilde{V}_h$. By coercivity and the fact that $a(\widetilde{u}_h;\widetilde{u}_h,\widetilde{v}_h)=a(u;u,\widetilde{v}_h)$ for every $\widetilde{v}_h \in \widetilde{V}_h$, and adding and subtracting $a(u; \widetilde{w}_h,\widetilde{e}_h)$,
\begin{equation}\label{eq:firsteq}
\begin{aligned}
\kappa_\text{min} |\widetilde{e}_h |^2_1 &\leq a(\widetilde{u}_h;\widetilde{e}_h,\widetilde{e}_h) = a(\widetilde{u}_h;\widetilde{w}_h, \widetilde{e}_h) - a(\widetilde{u}_h; \widetilde{u}_h,\widetilde{e}_h) =  I_1 + I_2,
\end{aligned}
\end{equation}
where
$$
I_1 = a(\widetilde{u}_h;\widetilde{w}_h, \widetilde{e}_h) - a(u;\widetilde{w}_h, \widetilde{e}_h)
\text{ and }
I_2 = a(u;\widetilde{w}_h-u, \widetilde{e}_h).
$$

Because $\kappa_j$ is Lipschitz continuous and using Cauchy-Schwarz inequality and H\"{o}lder inequality,
\begin{equation} \label{eq:2ndterm}
\begin{aligned}
I_1 &\leq \sum_{j=0}^{m_\Gamma} \int_{\Omega_j} |\kappa_j(\widetilde{u}_h) - \kappa_j(u)| \, |\widetilde{w}_h^{\, \prime} | \, |\widetilde{e}_h^{\, \prime} | \, \d x\\
&\leq \sum_{j=0}^{m_\Gamma} \int_{\Omega_j}C_{L,j} |u - \widetilde{u}_h| \, |\widetilde{w}^{\, \prime}_h| \, |\widetilde{e}^{\, \prime}_h|\d x\\
&\leq C_L \int_{\Omega} |(u - \widetilde{u}_h) \widetilde{w}_h^{\, \prime} | \, |\widetilde{e}_h^{\, \prime}|\d x\\
&\leq C_L \| (u - \widetilde{u}_h) \widetilde{w}_h^{\, \prime} \| \, | \widetilde{e}_h |_1\\
&\le C_L \| u -  \widetilde{u}_h \|_{0,3} \| \,  \widetilde{w}_h^{\, \prime}  \|_{0,6} \, | \widetilde{e}_h |_1.
\end{aligned}
\end{equation} 
To proceed further, a bound for $\| u - \widetilde{u}_h \|_{L^3(\Omega)}$ is desired. To simplify the presentation, set $e = u - \widetilde{u}_h$. By Cauchy-Schwarz inequality and H\"{o}lder inequality,
\begin{equation} \label{eq:CCC}
\|e\|_{L^3(\Omega)}^3 \leq \|e\| \|e^2\| \leq \|e\| \, \|e\|_{0,3}\|e\|_{0,6}.
\end{equation}
By Sobolev embedding theorem  (see for example p. 85 of \cite{200379}),  $H^1(\Omega) \hookrightarrow L^6(\Omega)$, which implies $\|e\|_{L^6(\Omega)} \leq C_{\rm{sl}} |e|_1$. With this, the inequality in \cref{eq:CCC} yields
\begin{equation*}
\|e\|_{L^3(\Omega)}^2 \leq C_{\rm{sl}} \|e\| \,  |e|_1. 
\end{equation*}
Using this in \cref{eq:2ndterm} gives
\begin{equation*}
I_1 \le C_L \sqrt{C_{\rm{sl}}  \|u -  \widetilde{u}_h\| \,  |u -  \widetilde{u}_h|_1} \,  \| \widetilde{w}_h^{\, \prime}  \|_{0,6} \, | \widetilde{e}_h |_1.
\end{equation*}
Notice also, due to the boundedness of $a(\cdot ; \cdot,\cdot)$,
\begin{equation*}
 I_2  \le \kappa_\text{max} | u -\widetilde{w}_h |_1 \,  | \widetilde{e}_h |_1.
\end{equation*}
Using all these estimates in \cref{eq:firsteq} gives
\begin{equation*}
\kappa_\text{min} |\widetilde{e}_h |_1 \le C_L \sqrt{C_{\rm{sl}}  \|u -  \widetilde{u}_h\| \,  |u -  \widetilde{u}_h|_1} \,  \| \widetilde{w}_h^{\, \prime}  \|_{0,6}  + \kappa_\text{max} | u -\widetilde{w}_h |_1,
\end{equation*}
from which we obtain
\begin{equation*}
|\widetilde{e}_h |_1 \le C \Big( \sqrt{  \| \widetilde{w}_h^{\, \prime}  \|^2_{0,6}  \|u -  \widetilde{u}_h\| \,  |u -  \widetilde{u}_h|_1} + | u -\widetilde{w}_h |_1 \Big), 
\end{equation*}
where
$$
C = \kappa_{\text{min}}^{-1} \max \big( C_L C_{\text{sl}}^{1/2}, \kappa_{\text{max}} \big).
$$
Putting this last inequality back to \cref{eq:triangle} and applying inequality
$$
\sqrt{a b} \le \sqrt{(\delta^{-2} a^2 + \delta^2 b^2)/2} \le \frac{a}{\delta \sqrt{2}} + \frac{\delta b}{\sqrt{2}}, ~~a>0, b>0, \delta>0,
$$
gives
\begin{equation*}
\begin{aligned}
|u - \widetilde{u}_h|_1 &\leq (1+C) |u - \widetilde{w}_h|_1 + C \sqrt{  \| \widetilde{w}_h^{\, \prime}  \|^2_{0,6}  \|u -  \widetilde{u}_h\| \,  |u -  \widetilde{u}_h|_1}\\
&\le (1+C) |u - \widetilde{w}_h|_1 + \frac{C}{\delta \sqrt{2}} \| \widetilde{w}_h^{\, \prime}  \|^2_{0,6}  \|u -  \widetilde{u}_h\| + \frac{\delta C}{\sqrt{2}}  |u -  \widetilde{u}_h|_1.
\end{aligned}
\end{equation*}
By choosing $\delta>0$ such that $\delta C/\sqrt{2} < 1$, estimate in \cref{eq:sebutsajadia} is established.
\end{proof}

Previous lemma quantifies approximation error $H^1$-seminorm in terms of approximation error in $L^2$-norm and the "quality" of $\widetilde{V}_h$. In the next lemma, the approximation error in $L^2$-norm is expressed in terms of approximation error in $H^1$-seminorm and yet another notion of quality of $\widetilde{V}_h$. The technique utilized to prove this lemma is a duality argument first introduced by Aubin-Nitche (see for example \cite{ciarlet2,braess_2007}). However, it was originally applied to variational formulations of linear boundary value problems. The duality argument relies on an adjoint problem associated with the aforementioned variational formulations.

Since \cref{eq:probstatment} and the associated variational formulation is nonlinear, a linearization is required that allows for construction of the adjoint problem that is linear. To this end, define $\mathcal{F} : H^1_0(\Omega) \to \mathbb{R}$ by $\mathcal{F}(v) = a(v; v,w)$ for every $w \in H^1_0(\Omega)$. The Fr\'{e}chet derivative of $\mathcal{F}$ at $v \in H^1_0(\Omega)$ is $\mathcal{F}^\prime(v) \in [H^1_0(\Omega)]^*$ such that
$$
\lim_{\eta \to 0 \in H^1_0(\Omega) } \| \mathcal{F}(v+\eta) - \mathcal{F}(v) - [\mathcal{F}^\prime(v)](\eta) \|_1 = 0.
$$
In this case,
$$
[\mathcal{F}^\prime(v)](\eta) = a(v ; \eta, w) + b(v; \eta, w), ~~\forall w \in H^1_0(\Omega),
$$
with
\begin{equation} \label{eq:myb}
b(v ; \eta, w) = \sum_{j=0}^{m_\Gamma} \int_{\Omega_j} D_2 \kappa_j(x,v) v^\prime(x) \eta(x) w^\prime(x)  \, {\rm d} x,
\end{equation}
where $D_2 \kappa_j$ is the partial derivative of $\kappa_j$ with respect to the second variable.
Now set $\sigma : [0,1] \to H^1_0(\Omega)$ by $\sigma(t) = \widetilde{u}_h + t (u-\widetilde{u}_h)$. By integral mean value theorem, 
\begin{equation} \label{eq:mvt}
\begin{aligned}
a(u; u,w) - a(\widetilde{u}_h;\widetilde{u}_h,w) &=
\mathcal{F}(u) - \mathcal{F}(u_h) \\
&= \int_0^1 [\mathcal{F}^\prime(\sigma(t))] (u - \widetilde{u}_h) \,\d t \\
&=\int_0^1 \big( a(\sigma(t) ; u - \widetilde{u}_h,w) + b(\sigma(t) ;  u - \widetilde{u}_h,w)  \big) \,\d t \\
&= \overline{a}(\sigma ; u - \widetilde{u}_h,w) + \overline{b}(\sigma ;  u - \widetilde{u}_h,w),  
\end{aligned}
\end{equation}
where
\begin{equation*}
\begin{aligned}
\overline{a}(\sigma ; v,w) &= \sum_{j=0}^{m_\Gamma} \int_{\Omega_j} \Bigg( \int_0^1 \kappa_j(x, [\sigma(t)](x)) \, {\rm d} t  \Bigg)v^\prime(x) w^\prime(x) \, {\rm d} x,\\
\overline{b}(\sigma ; v, w) &= \sum_{j=0}^{m_\Gamma} \int_{\Omega_j} \Bigg( \int_0^1 D_2 \kappa_j(x,[\sigma(t)](x)) [\sigma(t)]^\prime(x) \, {\rm d} t \Bigg) v(x) w^\prime(x)  \, {\rm d} x.
\end{aligned}
\end{equation*}
The above forms are linear in the second and third argument so both of them are bilinear form.
Given $\psi \in L^2(\Omega)$, the adjoint problem is to seek $\varphi \in H^1_0(\Omega)$ that is governed by
\begin{equation}  \label{eq:weakadjoint} 
\overline{a}(\sigma ; v,\varphi) +  \overline{b}(\sigma ; v, \varphi) = (v, \psi), ~~\forall v \in H^1_0(\Omega).
\end{equation}
Notice that this is a variational formulation of a linear boundary value problem. Under an additional assumption that $\kappa_j \in C^1(\overline{\Omega}_j \times \mathbb{R})$ and using standard tools from ordinary differential equations, existence of such a $\varphi$ is established in the Appendix.
\begin{lemma} \label{L2ineq}
Assume further that $\kappa_j \in C^1(\overline{\Omega}_j \times \mathbb{R})$ for every $j=1,\cdots,m_\Gamma$. There exists a positive contant $C>0$ independent of $u$ and $h$ such that
\begin{equation} \label{eq:mawar}
\|u - \widetilde{u}_h\|^2 \leq C | \varphi - \widetilde{w}_h |_1 \, |u - \widetilde{u}_h|_1, ~~\forall \widetilde{w}_h \in \widetilde{V}_h,
\end{equation}
where $\varphi \in H^1_0(\Omega)$ satisfies \cref{eq:weakadjoint}.
\end{lemma}
\begin{proof}
Let $e = u - \widetilde{u}_h$ and use \cref{eq:weakadjoint} with $\psi = e$ and $v = e$ and \cref{eq:mvt}
 to get
\begin{equation*}
\begin{aligned}
\| e \|^2 = (e,e) &= \overline{a}(\sigma ; e,\varphi) +  \overline{b}(\sigma ; e, \varphi)  = a(u; u,\varphi) - a(\widetilde{u}_h;\widetilde{u}_h,\varphi).
\end{aligned}
\end{equation*}
Using $a(\widetilde{u}_h;\widetilde{u}_h,\widetilde{w}_h) - a(u; u,\widetilde{w}_h) = 0$
for any $\widetilde{w}_h \in \widetilde{V}_h$ and add and subtract $a(u ; \widetilde{u}_h,\varphi - \widetilde{w}_h)$,
\begin{equation} \label{eq:DDD}
\begin{aligned}
\|e\|^2 = a(u ; u, \varphi - \widetilde{w}_h) - a(\widetilde{u}_h ; \widetilde{u}_h, \varphi - \widetilde{w}_h) = J_1 + J_2,
\end{aligned}
\end{equation}
where
$$
J_1 = a(u ; u-\widetilde{u}_h, \varphi - \widetilde{w}_h),  \text{ and }
J_2 = a(u ; \widetilde{u}_h, \varphi - \widetilde{w}_h)  - a(\widetilde{u}_h ; \widetilde{u}_h, \varphi - \widetilde{w}_h). 
$$
Using the boundedness of $a(\cdot ; \cdot,\cdot)$,
\begin{equation*}
J_1 \le \kappa_{\text{max}} |u-\widetilde{u}_h|_1 \, |\varphi - \widetilde{w}_h|_1.
\end{equation*}
By applying the Lipschitz continuity of $\kappa_j$, Cauchy-Schwarz inequality, the boundedness of $\widetilde{u}_h$, i.e., $\|\widetilde{u}_h\|_1 \leq C_0^{-1} \| f \|$, and  $\|w\|_{0,\infty} \leq C_\Omega |w|_1$
for any $w \in H^1_0(\Omega)$, $J_2$ is estimated as follows:
\begin{equation*}
\begin{aligned}
J_2 &\le \sum_{j=0}^{m_\Gamma} \int_{\Omega_j} |\kappa_j(u) - \kappa_j(\widetilde{u}_h)| \, |u_h^\prime | \, |(\varphi - \widetilde{w}_h)^\prime| \, \d x \\
&\le \sum_{j=0}^{m_\Gamma} \int_{\Omega_j} C_{L,j} | u - \widetilde{u}_h| \, |\widetilde{u}_h^\prime | \, |(\varphi - \widetilde{w}_h)^\prime| \, \d x \\
&\leq C_L \, \| u- \widetilde{u}_h \|_{0,\infty} \, |\widetilde{u}_h|_1 \, |\varphi - \widetilde{w}_h |_1\\
&\leq C_L  C_0^{-1} C_\Omega \| f \| \, | u- \widetilde{u}_h |_1 \, |\varphi - \widetilde{w}_h |_1.
\end{aligned}
\end{equation*}
Putting these estimates back into \cref{eq:DDD} gives the desired result.
\end{proof}

At this stage, the tools needed to derive an error estimate of the finite element approximations are in place. Based on the results in \Cref{H1ineq} and \Cref{L2ineq}, an issue here is the quantification of $| u - \widetilde{w}_h |_1$, $| \widetilde{w}_h |_{1,6}$, and $| \varphi - \widetilde{w}_h |_1$, so it boils down to the approximation properties of the finite element spaces $\widetilde{V}_h = V_h^p$ or $\widetilde{V}_h = V_{h,\tE}^p$. In fact, this is where $V_h^p$ behaves differently from  $V_{h, \tE}^p$, in a sense that the interpolation operator $\mathcal{I}^p_h$ for $V^p_h$ does not have analog approximation properties of $\mathcal{I}^p_{h,\tE}$ for $V^p_{h,\tE}$ as described in \Cref{global_int} and \Cref{lem:norm16}. The next theorem is the error estimate for SGFEM solution.

\begin{theorem}\label{errorestimates}(Error estimates)
Assume further that $\kappa_j \in C^1(\overline{\Omega}_j \times \mathbb{R})$ for every $j=1,\cdots,m_\Gamma$. If $u \in H^{p+1}_\Gamma(\Omega)$ and $u_{h, \emph \tE} \in V^p_{h,\emph \tE}$ is its SGFEM approximation, then there exists an $h_0>0$ such that for any $h<h_0$, 
\begin{equation} \label{eq:errorestimate}
\|u - u_{h, \emph \tE}\| + h |u - u_{h,\emph \tE}|_1 \leq Ch^{p+1} \Bigg(\sum_{j=0}^{m_\Gamma} |u|^2_{p+1,\Omega_j}\Bigg)^{\frac{1}{2}},
\end{equation}
for some $C>0$.
\end{theorem}
\begin{proof}
Here we mainly employ \Cref{H1ineq} and \Cref{L2ineq} with $\widetilde{u}_h = u_{h,\tE}, \widetilde{w}_h = w_{h,\tE}$ all belonging to $\widetilde{V}_h = V^p_{h,\tE}$. Thus, by \Cref{H1ineq},
\begin{equation} \label{eq:1sterrE}
|u - u_{h,\tE}|_1 \leq C( |u - w_{h,\tE} |_1 + | w_{h,\tE} |^2_{1,6} \|u - u_{h,\tE}\|) ~\forall ~w_{h,\tE} \in V^p_{h,\tE}.
\end{equation}

By choosing $w_{h,\tE} = \mathcal{I}^p_{h,\tE} u$ and using \Cref{lem:norm16},
$$
| \mathcal{I}^p_{h,\tE} u |_{1,6} \le |u|_{1,6} + | u - \mathcal{I}^p_{h,\tE} u |_{1,6} \le 
|u|_{1,6} + h^{p-\frac{1}{3}}  \Bigg(\sum_{j=0}^{m_\Gamma}|u|^2_{p+1,\Omega_j}\Bigg)^{\frac{1}{2}}.
$$
Moreover, by Sobolev embedding theorem (see for example p. 85 of \cite{200379}),
$$
|u|_{1,6} = \Bigg( \sum_{j=0}^{m_\Gamma} |u|^6_{1,6,\Omega_j} \Bigg)^{\frac{1}{6}}
\le  \Bigg( \sum_{j=0}^{m_\Gamma} C_j ||u||^6_{p+1,\Omega_j} \Bigg)^{\frac{1}{6}} \le
C \Bigg( \sum_{j=0}^{m_\Gamma}  ||u||^2_{p+1,\Omega_j} \Bigg)^{\frac{1}{2}}.
$$
Using all these estimates and \Cref{global_int} in \Cref{eq:1sterrE} give
\begin{equation} \label{eq:h1errnow}
|u - u_{h,\tE}|_1 \leq Ch^p \Bigg(\sum_{j=0}^{m_\Gamma}|u|^2_{p+1,\Omega_j}\Bigg)^{\frac{1}{2}} + C \Bigg( \sum_{j=0}^{m_\Gamma}  ||u||^2_{p+1,\Omega_j} \Bigg)^{\frac{1}{2}} \|u - u_{h,\tE}\|.
\end{equation}

Next, let $\varphi \in H^1_0(\Omega)$ be the solution of \cref{eq:weakadjoint} (see \Cref{pr:A1} for its existence). By choosing $\widetilde{w}_h = w_{h,\tE} = \mathcal{I}^1_{h,\tE} \varphi \in V^p_{h,\tE}$ in \Cref{L2ineq} and using \Cref{pr:A2}, we obtain
\begin{equation} \label{eq:h1l2}
\|u - u_{h,\tE} \| \leq C h \, |u - u_{h,\tE}|_1.
\end{equation}
Fix a sufficiently small $h_0>0$ such that 
$$
C h_0 \Bigg( \sum_{j=0}^{m_\Gamma}  ||u||^2_{p+1,\Omega_j} \Bigg)^{\frac{1}{2}} < 1.
$$
Now with $h < h_0$, we may put  \cref{eq:h1l2} into \cref{eq:h1errnow} and combine the last term on the right hand side with the term on the left hand side to give
$$
|u - u_{h,\tE}|_1 \leq Ch^p \Bigg(\sum_{j=0}^{m_\Gamma}|u|^2_{p+1,\Omega_j}\Bigg)^{\frac{1}{2}}.
$$
This last inequality and \cref{eq:h1l2} yield the desired result.
\end{proof}
\section{Numerical Examples}
\setcounter{equation}{0}
\label{sec:NE}
In this section we give detailed numerical examples to demonstrate the application of SGFEM approximation to some quasilinear elliptic problems whose analytic solutions are available.  We also aim to show that the approximation errors in the numerical examples reflect the optimal convergence properties of SGFEM as established in \Cref{errorestimates}. To solve the variational formulation \cref{eq:finitevariational}, standard Newton's method of iteration is employed:
\begin{algorithm}
\caption{} \label{alg:1}
\begin{algorithmic}
\State Set $\widetilde{u}_h^{(0)} \in \widetilde{V}_h$ (an initial guess).
\For{$n=1,2,\cdots, \text{until convergence}$}
\State Find $\delta_h \in \widetilde{V}_h$ governed by
\begin{equation}\label{eq:femformulation}
a(\widetilde{u}_h^{(n-1)}  ; \delta_h,\widetilde{w}_h) + b(\widetilde{u}_h^{(n-1)} ; \delta_h,\widetilde{w}_h) = \ell(\widetilde{w}_h) - a(\widetilde{u}_h^{(n-1)};\widetilde{u}_h^{(n-1)},\widetilde{w}_h), \hspace{0.5cm} \widetilde{w}_h \in \widetilde{V}_h.
\end{equation}
\State Set $\widetilde{u}_h^{(n)}  = \widetilde{u}_h^{(n-1)} + \delta_h$.
\EndFor
\end{algorithmic}
\end{algorithm}
The form $b(\cdot ; \cdot,\cdot)$ in \cref{eq:femformulation} is as expressed in \cref{eq:myb}.
 For all examples below, the initial guess is $\widetilde{u}^{(0)}_h = 0$. Denoting the linear system associated with \cref{eq:femformulation} by $\mathcal{A}^{(n-1)} \boldsymbol{\delta} = \boldsymbol{r}^{(n-1)}$,
a convergence is declared when $\| \boldsymbol{r}^{(n-1)}\|_{\infty} \leq 10^{-10}$.

\subsection{Example 1 (A Quasilinear problem with 2 interfaces)}
\begin{figure}[H] 
\centering
\includegraphics[scale = 0.65]{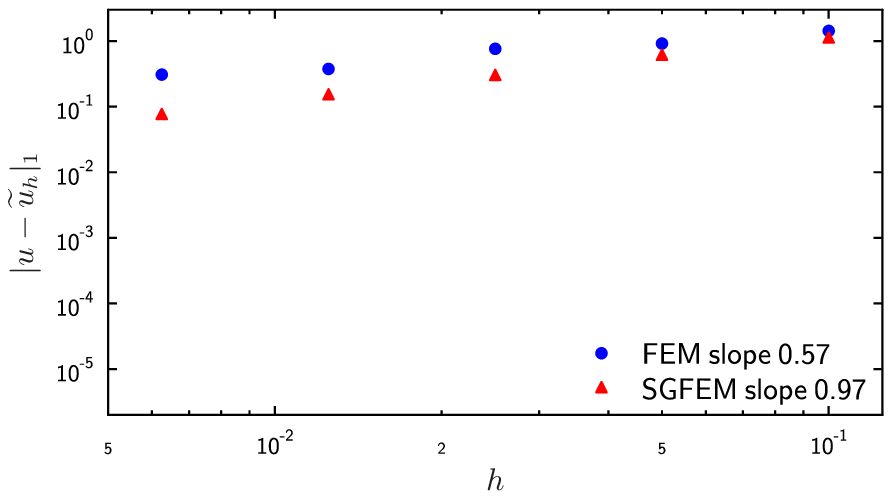}
\includegraphics[scale = 0.65]{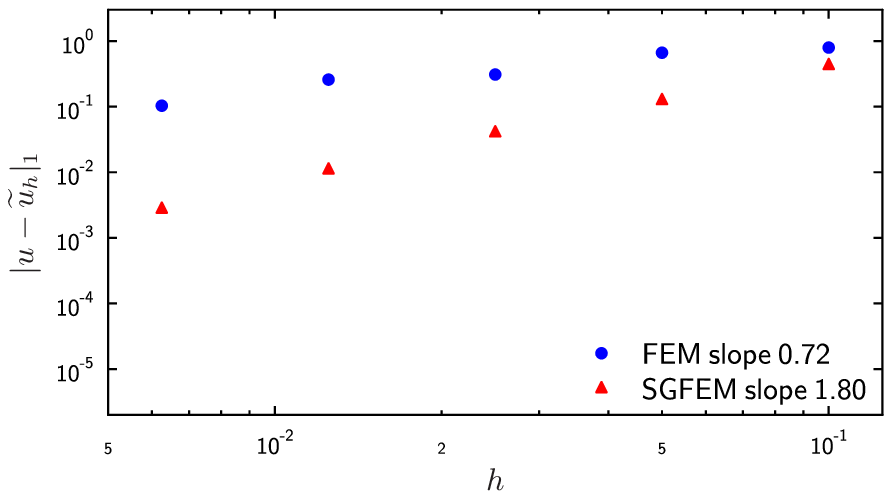}\\
\includegraphics[scale = 0.65]{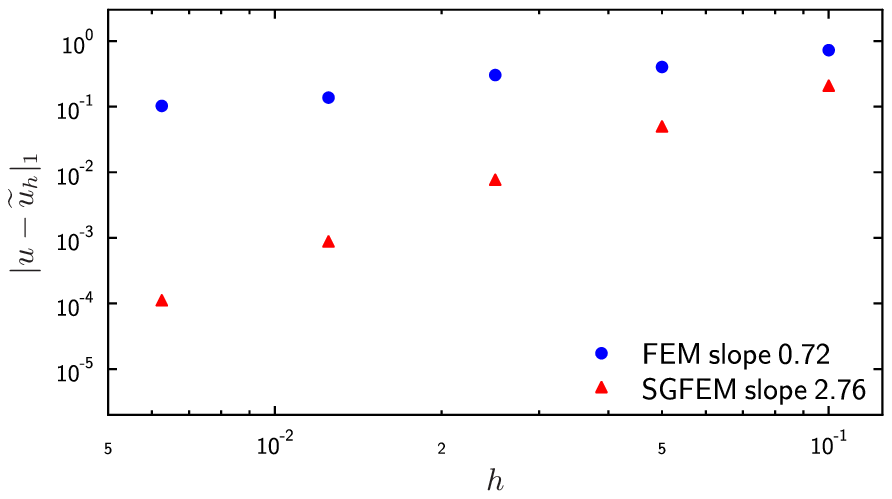}
\includegraphics[scale = 0.65]{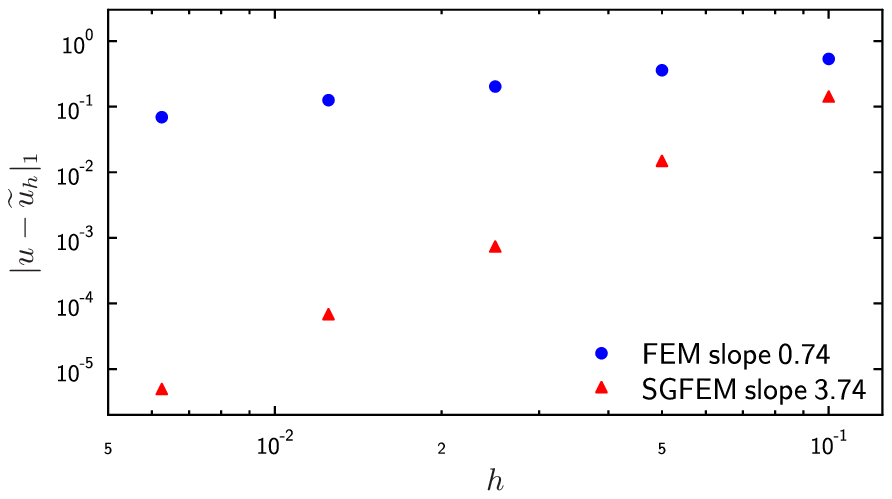}
\caption{Example 1: $|u - \widetilde{u}_h|_1$ v.s. $h$ using $p=1$ (top left), $p=2$ (top right),  $p=3$ (bottom left), and $p=4$ (bottom right). FEM corresponds to $\widetilde{u}_h = u_h \in V_h^p$ and SGFEM corresponds to $\widetilde{u}_h = u_{h,\tE} \in V_{h,\tE}^p$.}
\label{figtwo}
\end{figure}

\begin{figure}[H]
\centering
\includegraphics[scale =0.65]{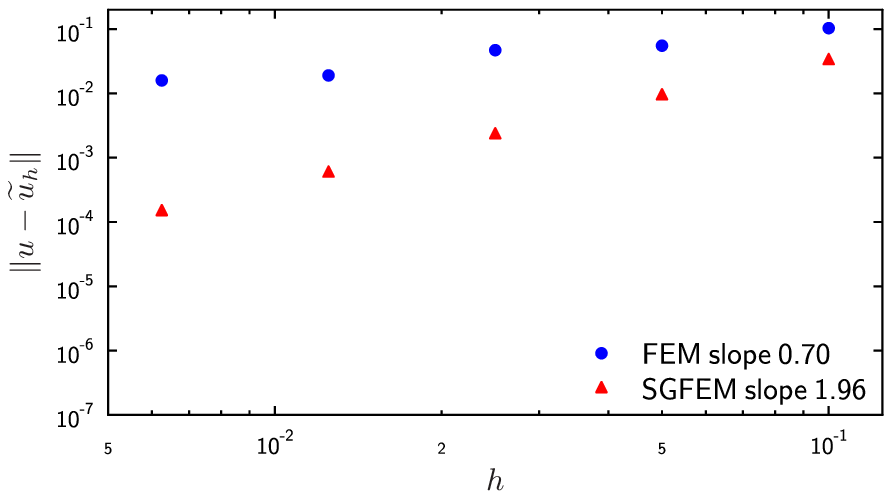}
\includegraphics[scale =0.65]{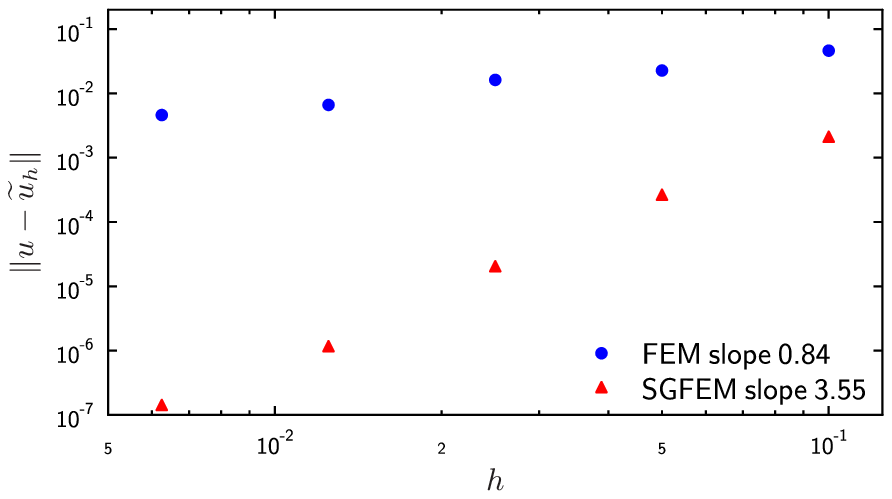}
\caption{Example 1: $\|u - \widetilde{u}_h\|$ v.s. $h$ using $p=1$ (left) and $p=3$ (right). FEM corresponds to $\widetilde{u}_h = u_h \in V_h^p$ and SGFEM corresponds to $\widetilde{u}_h = u_{h,\tE} \in V_{h,\tE}^p$.}
\label{figL2norm21}
\end{figure}

Let $\Omega = (0,1)$ and the interfaces $\gamma_1 = 1/3$, $\gamma_2 = 2/3$. We choose $f(x) = 5x$, and
\begin{equation} \label{eq:prob1}
\kappa(x,u) = 
\begin{cases}
e^{a_0u}, &x\in(0, \gamma_1],\\
e^{a_1u}, & x \in (\gamma_1, \gamma_2], \\
e^{a_2u}, & x \in(\gamma_2,1). 
\end{cases}
\end{equation} 
The solution is obtained by applying fundamental theorem of calculus and imposing the continuity of the solution and the flux. The analytical solution is given in \cref{eq:anasol1}.
\begin{figure}[H] 
\centering
\includegraphics[scale=0.7]{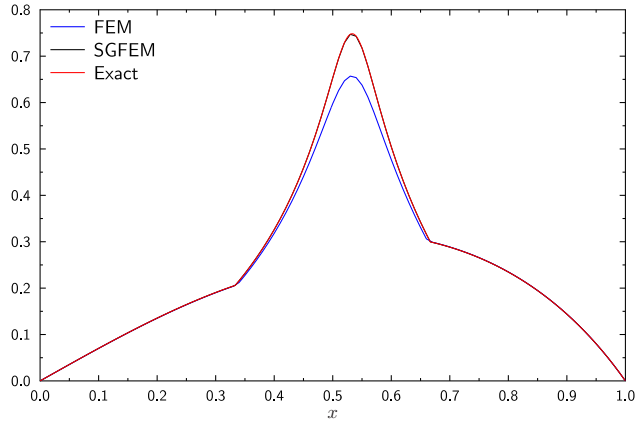}
\caption{Example 1: A comparison of FEM and SGFEM using $N=100$ and $p=1$.} 
\label{figone}
\end{figure}
Profile plots of the true solution and the approximate solutions using linear FEM and SGFEM are given in \Cref{figone}, with $a_0  = 0.01$, $a_1 = -6$, and $a_2 = 1$.
This shows a better performance of SGFEM to approximate the interface problem compared to the standard FEM. The contrast ratio for this example, which is calculated as $\kappa_{\text{max}}/\kappa_{\text{min}}$ is equal to 120. A uniform discretization of the domain into $N = 10, 20, 40, 80, 160$ elements is done in such way that the resulting mesh configurations guarantees that an interface $\gamma \in \Gamma$ is always located inside an element $\tau \in \mathcal{T}_{h,\Gamma}$. Comparison of the errors is shown in log-log  plots in \Cref{figtwo} for $H^1$ semi-norm and in  \Cref{figL2norm21} for $L^2$ norm.  The slopes are given in the plots for each approximation to see the convergence rate. Plots in these figures confirm the optimality of convergence property of SGFEM solution.

\subsection{Example 2 (Quasilinear problem with 3 interfaces)} 
\begin{figure}[H] 
\centering
\includegraphics[scale = 0.7]{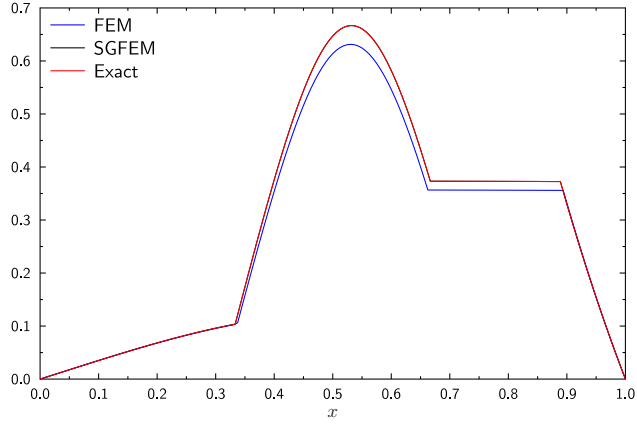}
\caption{Example 2: A comparison of FEM and SGFEM using $N=100$ and $p=1$.}
\label{figthree}
\end{figure}

\begin{figure}[H] 
\centering
\includegraphics[scale =0.65]{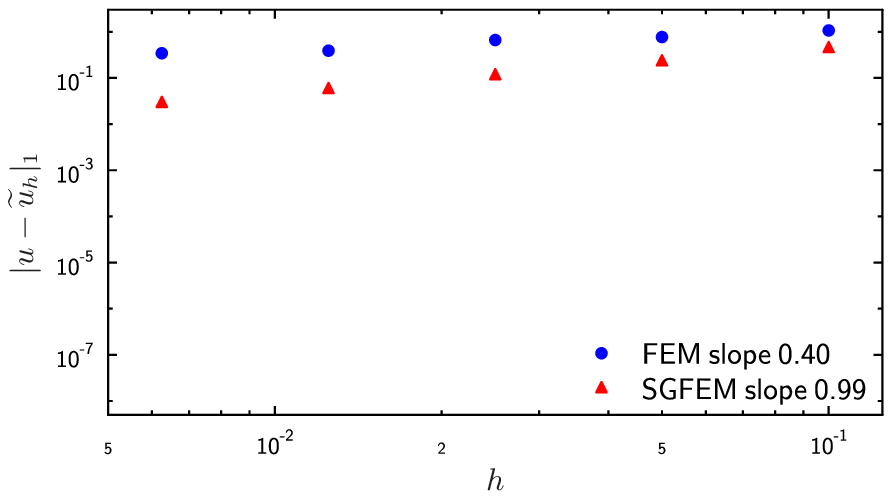}
\includegraphics[scale =0.65]{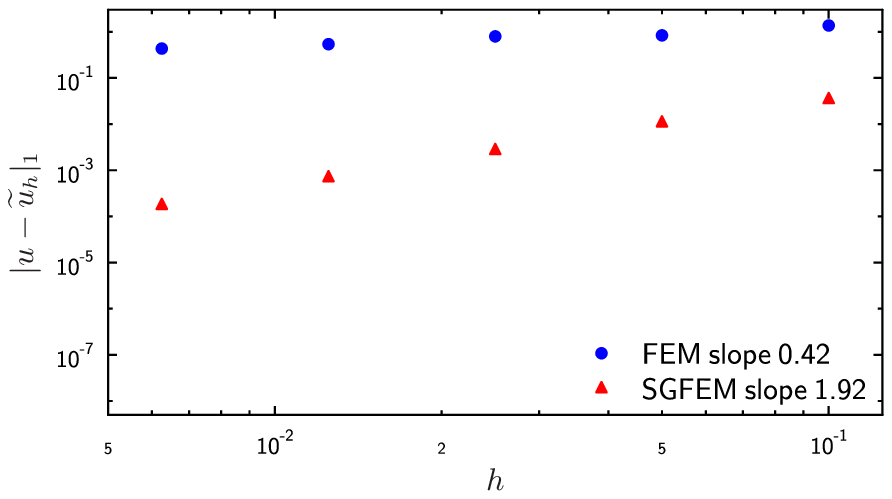}\\
\includegraphics[scale =0.65]{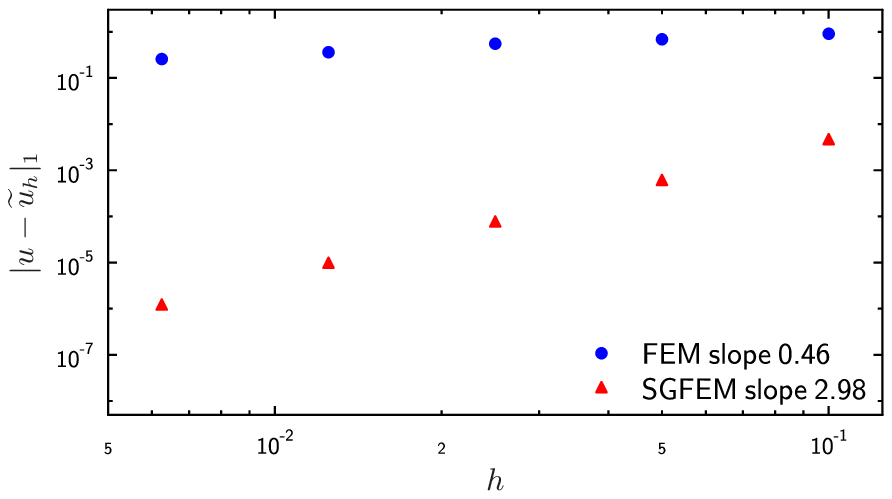}
\includegraphics[scale =0.65]{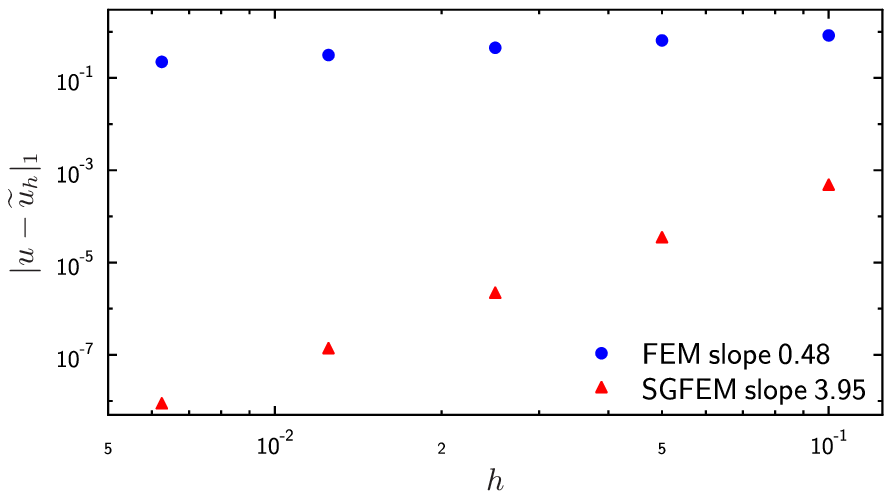}
\caption{Example 2: $|u - \widetilde{u}_h|_1$ v.s. $h$ using $p=1$ (top left), $p=2$ (top right),  $p=3$ (bottom left), and $p=4$ (bottom right). FEM corresponds to $\widetilde{u}_h = u_h \in V_h^p$ and SGFEM corresponds to $\widetilde{u}_h = u_{h,\tE} \in V_{h,\tE}^p$.}
\label{figfour}
\end{figure}

\begin{figure}[H]
\centering
\includegraphics[scale = 0.65]{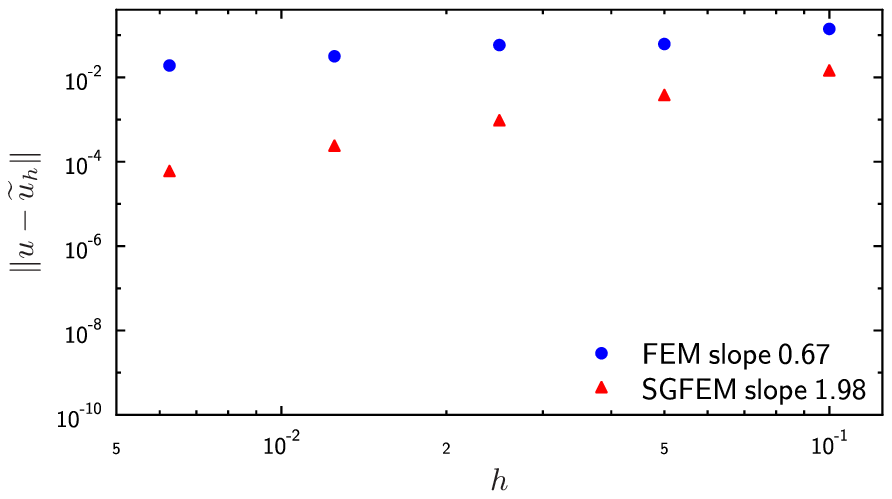}
\includegraphics[scale = 0.65]{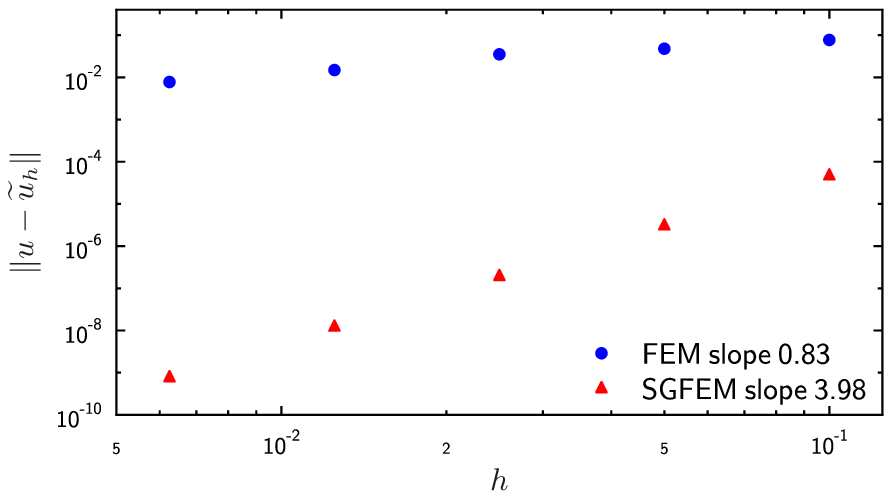}
\caption{Example 2: $\|u - \widetilde{u}_h\|$ v.s. $h$ using $p=1$ (left) and $p=3$ (right). FEM corresponds to $\widetilde{u}_h = u_h \in V_h^p$ and SGFEM corresponds to $\widetilde{u}_h = u_{h,\tE} \in V_{h,\tE}^p$.}
\label{figL2norm16}
\end{figure}

Let $\Omega = (0,1)$ and the interfaces $\gamma_1 = 1/3, \gamma_2 = 2/3$, and $\gamma_3 = 8/9$. We choose $f(x) = \sin(\pi x)$ and 
\begin{equation}  \label{eq:prob2}
\kappa(x,u) = 
\begin{cases}
a_0e^{-u}, & x \in (0,\gamma_1],\\
a_1e^{-u}, &x\in(\gamma_1,\gamma_2],\\
a_2e^{-u}, & x \in(\gamma_2,\gamma_3],\\
a_3e^{-u}, &x\in(\gamma_3,1).
\end{cases}
\end{equation}
The true solution is given in \cref{eq:anasol2}. Similar to Example 1, we plot the true solution and the approximate solutions of linear FEM and SGFEM in \Cref{figthree} using $a_0 = 1$, $a_1 = 0.05$, $a_2 =100$, and $a_3 = 0.1$. The approximate solutions are produced under the same mesh configurations as Example 1. 
The contrast ratio for this example is 2684. The corresponding errors in $H^1$ semi-norm and in $L^2$ norm are respectively plotted in \Cref{figfour} and \Cref{figL2norm16}. Again, results in this example validates the optimal convergence property of SGFEM. It is also observed that as the contrast coefficients become higher, the optimal convergence for SGFEM with $p=1$ is still preserved, however the convergence rate of SGFEMs for $p>1$ is deteriorating. 


\section{FEMs with Local Conservation Constraints}
\setcounter{equation}{0}
\label{sec:FELC}
In this section, we devise a procedure to construct an approximate solution that in addition to satisfying \cref{eq:finitevariational}, it also obeys a mass balance (commonly called local conservation) over a set of control volumes of $\Omega$. A collection of $N^*$ control volumes is denoted by $\mathcal{T}^* = \{ (t_{i-1},t_i) \subset \overline{\Omega} : i=1,\cdots,N^* \}$. We assume that $N^* < \text{dim}(\widetilde{V}_h)$, where as before either $\widetilde{V}_h=V_h^p$ or  $\widetilde{V}_h =V^p_{h,\tE}$.
Several relevant examples of the control volumes are
\begin{equation} \label{eq:CV}
\begin{aligned}
\mathcal{T}^* &= \{ \Omega \},\\
\mathcal{T}^* &= \{ \Omega_j : j=0,\cdots,m_\Gamma \}, \text{ and }\\
\mathcal{T}^* &= \{ (t_j,t_{j+1}) \subset \Omega, j=1,\cdots,N-1 : t_i = \text{midpoint of } \tau_i \in\mathcal{T}_h, i=1,\cdots,N \}.
\end{aligned}
\end{equation}

Now let $\kappa : \Omega \times \mathbb{R} \to \mathbb{R}$ be defined such that $\kappa |_{\Omega_j} = \kappa_j$ for $j =1,\cdots, m_\Gamma$. For a $\tau^* = (t_l,t_r) \in \mathcal{T}^*$, define
\begin{equation*}
\begin{aligned}
C_{\tau^*}: H^1_0(\Omega) \times H^1_0(\Omega) \to \mathbb{R} \text{ as } C_{\tau^*}(v;w) &= -\kappa(x, v(x))  w^\prime(x) \Big |_{t_l}^{t_r}  \text{ and }\\
\ell_{\tau^*} &=  \int_{t_l}^{t_r} f(x) \, \d x.
\end{aligned}
\end{equation*}
Notice that $C_{\tau^*}$ is linear with respect to the second argument. Any approximate solution $\widetilde{u}$ is locally conservative if
\begin{equation} \label{eq:conservation}
C_{\tau^*}( \widetilde{u} ; \widetilde{u} ) = \ell_{\tau^*} \text{ for every } \tau^* \in \mathcal{T}^*.
\end{equation}
This property is obviously not satisfied by $\widetilde{u}_h \in \widetilde{V}_h$ governed by \cref{eq:finitevariational}.  Formally, the intention is to seek $\widetilde{u}_h \in \widetilde{V}_h$ governed by
  \begin{equation*}
  \begin{cases}
  a(\widetilde{u}_h ; \widetilde{u}_h, \widetilde{w}_h) = \ell(\widetilde{w}_h), ~~\forall \widetilde{w}_h \in  \widetilde{V}_h,\\
  C_{\tau^*}( \widetilde{u}_h ; \widetilde{u}_h ) = \ell_{\tau^*} \text{ for every } \tau^* \in \mathcal{T}^*.
  \end{cases}
  \end{equation*}
 Unfortunately, posing it this way is practically infeasible because the above system has more equations to satisfy than the number of unknowns involved in it.
 
 To get a way around this obstacle, we adopt a Lagrange multiplier technique introduced in \cite{Presho15,Abreu17}. The main idea lies on a recognition that typical linear variational formulation problem posed on a Banach space is equivalent to a minimization of a certain energy functional over that same space. By introducing a set of Lagrange multipliers, a new functional is created to include a set of constraints. The problem becomes a minimization of this new functional.
 
 However, unlike the usual linear variational formulation, the nonlinearity in $a(v ; v, w)$ does not allow for a direct energy functional that can be minimized. To tackle this issue, we propose to perform a minimization of a linear functional that results from a linearized variational formulation, which is then equipped with the Lagrange multipliers to include the constraints. We note that as described in \Cref{sec:NE},
a linearization based on the Fr\'{e}chet derivative of $a(v; v, w)$ has been enforced to allow for the implementation of Newton's method of iteration.

To describe the conceptual framework, fix $z \in H^1_0(\Omega)$ and set
\begin{equation*}
\mathcal{J}(z ; w, \boldsymbol{\zeta}) = \frac{1}{2} a(z ; w,w) - \ell(w) + \sum_{\tau^* \in \mathcal{T}^*} \zeta_{\tau^*} (C_{\tau^*}(z ; w) - \ell_{\tau^*}),
\end{equation*}
for any $w \in H^1_0(\Omega)$ and $\boldsymbol{\zeta} \in \mathbb{R}^{N^*}$ whose components are $\zeta_{\tau^*} \in \mathbb{R}$.  The Fr\'{e}chet derivative of $\mathcal{J}$ at $(v,\boldsymbol{\lambda})$ is denoted by a bilinear functional $\mathcal{J}^\prime(z ; v,\boldsymbol{\lambda}) : H^1_0(\Omega) \times \mathbb{R}^{N^*} \to \mathbb{R}$ and is defined as
\begin{equation*}
[\mathcal{J}^\prime(z ; v,\boldsymbol{\lambda})](w,\boldsymbol{\zeta}) = a(z ; v,w) - \ell(w) +
\sum_{\tau^* \in \mathcal{T}^*} \lambda_{\tau^*} C_{\tau^*}(z ; w) +
\sum_{\tau^* \in \mathcal{T}^*} \zeta_{\tau^*} (C_{\tau^*}(z ; v) - \ell_{\tau^*}).
\end{equation*}
%
%
If there is $(\widetilde{u}_h, \boldsymbol{\lambda}) \in \widetilde{V}_h \times \mathbb{R}^{N^*}$  such that
\begin{equation} \label{eq:Jp0}
[\mathcal{J}^\prime(z ; \widetilde{u}_h,\boldsymbol{\lambda})](\widetilde{w}_h,\boldsymbol{\zeta}) = 0, ~\forall (\widetilde{w}_h, \boldsymbol{\zeta}) \in \widetilde{V}_h \times \mathbb{R}^{N^*},
\end{equation}
then $\mathcal{J}(z ;\widetilde{u}_h, \boldsymbol{\lambda})$ is an extremal value. We end up with seeking $(\widetilde{u}_h, \boldsymbol{\lambda}) \in \widetilde{V}_h \times \mathbb{R}^{N^*}$ that is governed by \cref{eq:Jp0}, which is equivalent to
\begin{equation} \label{eq:finiteconstraintsform} 
\begin{aligned} 
&\text{ find }  (\widetilde{u}_h, \boldsymbol{\lambda}) \in \widetilde{V}_h \times \mathbb{R}^{N^*} \text{ that is governed by }\\
&\begin{cases}
a(z;  \widetilde{u}_h,\widetilde{w}_h) + \displaystyle \sum_{\tau^* \in \mathcal{T}^*} \lambda_{\tau^*} C_{\tau^*}(z ; \widetilde{w}_h))  = \ell(\widetilde{w}_h) , ~~ \forall \widetilde{w}_h \in \widetilde{V}_h\\
C_{\tau^*}(z ; \widetilde{u}_h) =  \ell_{\tau^*}, ~~\forall \tau^* \in \mathcal{T}^*.
\end{cases}
\end{aligned}
\end{equation}

The system  \cref{eq:finiteconstraintsform}  is placed within an iterative procedure written as follows:
\begin{algorithm}
\caption{} \label{alg:2}
\begin{algorithmic}
\State Set $\widetilde{u}_h^{(0)} \in \widetilde{V}_h$ (an initial guess).
\For{$n=1,2,\cdots, \text{until convergence}$}
\State Find $(\widetilde{u}_h^{(n)}, \boldsymbol{\lambda}^{(n)} ) \in \widetilde{V}_h \times \mathbb{R}^{N^*}$
governed by
\begin{equation*}
\begin{cases}
a(\widetilde{u}_h^{(n-1)}  ; \widetilde{u}_h^{(n)},\widetilde{w}_h) + \displaystyle \sum_{\tau^* \in \mathcal{T}^*} \lambda_{\tau^*}^{(n)}  C_{\tau^*}( \widetilde{u}_h^{(n-1)}  ; \widetilde{w}_h)  = \ell(\widetilde{w}_h), ~~ \widetilde{w}_h \in \widetilde{V}_h,\\
C_{\tau^*}(\widetilde{u}_h^{(n-1)} ; \widetilde{u}_h^{(n)}) =  \ell_{\tau^*}, ~~\forall \tau^* \in \mathcal{T}^*.
\end{cases}
\end{equation*}
\EndFor
\end{algorithmic}
\end{algorithm}

Notice that this algorithm is a fixed point type iteration. Supposing that the iteration converges to a limit $(\widetilde{u}_h, \boldsymbol{\lambda})$, namely $\| \widetilde{u}_h^{(n)} - \widetilde{u}_h \|_1 \to 0$ and $\| \boldsymbol{\lambda}^{(n)} - \boldsymbol{\lambda} \| \to 0$ as $n \to \infty$, and the limit satisfies
\begin{equation} \label{eq:CC}
\begin{cases}
a(\widetilde{u}_h  ; \widetilde{u}_h^,\widetilde{w}_h) + \displaystyle \sum_{\tau^* \in \mathcal{T}^*} \lambda_{\tau^*}  C_{\tau^*}( \widetilde{u}_h  ; \widetilde{w}_h)  = \ell(\widetilde{w}_h), ~~ \widetilde{w}_h \in \widetilde{V}_h,\\
C_{\tau^*}(\widetilde{u}_h ; \widetilde{u}_h) =  \ell_{\tau^*}, ~~\forall \tau^* \in \mathcal{T}^*,
\end{cases}
\end{equation}
then we may also apply a modified Newton's method of iteration to approximate $(\widetilde{u}_h, \boldsymbol{\lambda})$ in \cref{eq:CC}:

\begin{algorithm}
\caption{} \label{alg:3}
\begin{algorithmic}
\State Set $(\widetilde{u}_h^{(0)},\boldsymbol{\lambda}^{(0)}) \in \widetilde{V}_h \times \mathbb{R}^{N^*}$ (an initial guess).
\For{$n=1,2,\cdots, \text{until convergence}$}
\State Find $(\delta_h, \boldsymbol{\delta} ) \in \widetilde{V}_h \times \mathbb{R}^{N^*}$ governed by
\begin{equation*}
\begin{cases}
a(\widetilde{u}_h^{(n-1)}  ; \delta_h,\widetilde{w}_h) + \displaystyle \sum_{\tau^* \in \mathcal{T}^*} \delta_{\tau^*} [Q_{\tau^*}(\widetilde{u}_h^{(n-1)})](\widetilde{w}_h)  = R^{(n-1)}_1(\widetilde{w}_h)
, ~~ \widetilde{w}_h \in \widetilde{V}_h,\\
Q_{\tau^*}(\widetilde{u}_h^{(n-1)})](\delta_h) = R^{(n-1)}_2, ~~\forall \tau^* \in \mathcal{T}^*.
\end{cases}
\end{equation*}
\State Set $\widetilde{u}_h^{(n)} = \widetilde{u}_h^{(n-1)} + \delta_h$ and $\boldsymbol{\lambda}^{(n)} = \boldsymbol{\lambda}^{(n-1)} + \boldsymbol{\delta}$.
\EndFor
\end{algorithmic}
\end{algorithm}
Here,
\begin{equation*}
\begin{aligned}
[Q_{\tau^*}(z)](w) &= -\kappa(x, z(x))  w^\prime(x) - D_2 \kappa(x, z(x)) z^\prime(x) w(x) \Big |_{t_l}^{t_r}, ~~\forall w \in H^1_0(\Omega),\\
R^{(n-1)}_1(\widetilde{w}_h) &=  \ell(\widetilde{w}_h) - a(\widetilde{u}_h^{(n-1)}  ; \widetilde{u}_h^{(n-1)},\widetilde{w}_h), \\
R^{(n-1)}_2 &= \ell_{\tau^*} - C_{\tau^*}(\widetilde{u}_h^{(n-1)} ; \widetilde{u}_h^{(n-1)}).
\end{aligned}
\end{equation*}

\section{A Numerical Example for the Local Conservation} \label{sec:anelc}
\setcounter{equation}{0}
We use Example 2 in the previous section to compare the FEM/SGFEM solutions to the ones satisfying local conservation constraints to be imposed on a set control volumes $\mathcal{T}^*$ given by the third example in \cref{eq:CV}. We use \Cref{alg:3} to obtain $\widehat{u}_h \in V_h^p$ or $\widehat{u}_{h,\tE} \in V_{h,\tE}^p$ that satisfies the local conservation property for every $\tau^* \in \mathcal{T}^*$. The iteration was stopped once the relative residual was reduced by a factor of $10^{-10}$.

\begin{figure}[] 
\centering
\includegraphics[scale=0.65]{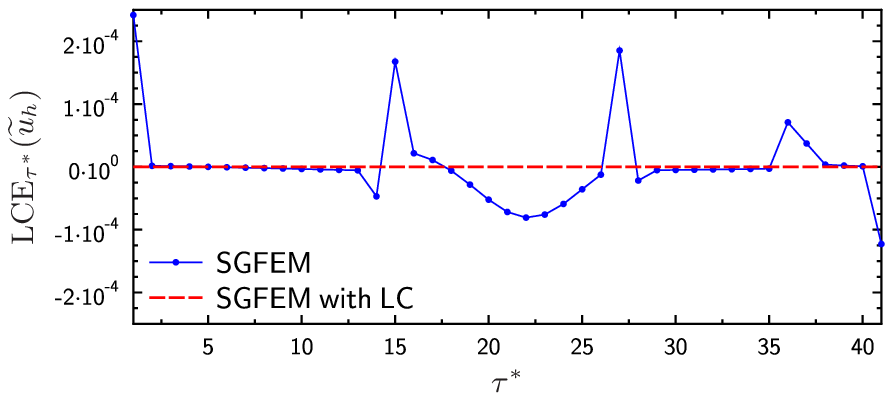}
\includegraphics[scale=0.675]{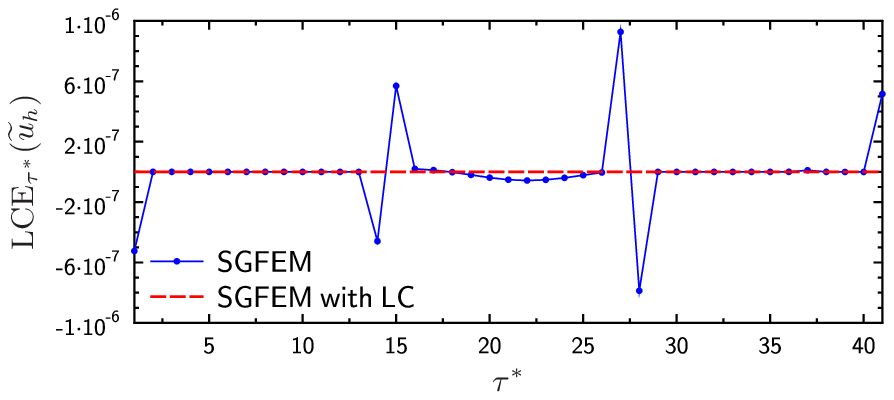}
\caption{Example 2: $\text{LCE}_{\tau^*}(\widetilde{u}_h)$ v.s. $\tau^*$ with $h =1/40$ and using $p=2$ (left) and $p=3$ (right). SGFEM corresponds to $\widetilde{u}_h = u_{h,\tE} \in V_{h,\tE}^p$, and SGFEM with LC corresponds to $\widetilde{u}_h = \widehat{u}_{h,\tE} \in V_{h,\tE}^p$ obtained from \Cref{alg:3}.}
\label{figlce}
\end{figure}

First, we calculate the local conservation errors in each $\tau^*$ to verify that the Lagrange multiplier technique indeed satisfies the conservation property. For this purpose we define the local conservation error (LCE) as
$$
\text{LCE}_{\tau^*} (w)= C_{\tau^*}(w;w) - \ell_{\tau^*} ,~~\tau^* \in \mathcal{T}^*.
$$
It is shown in \Cref{figlce} that $\mbox{LCE}_{\tau^*}(u_{h,\tE};u_{h,\tE}) \neq 0$, which confirms that the local conservation is violated, while $\mbox{LCE}_{\tau*}(\widehat{u}_{h,\tE};\widehat{u}_{h,\tE}) = 0$, which signifies that the local conservation is satisfied.

\begin{figure}[]
\centering
\includegraphics[scale=0.65]{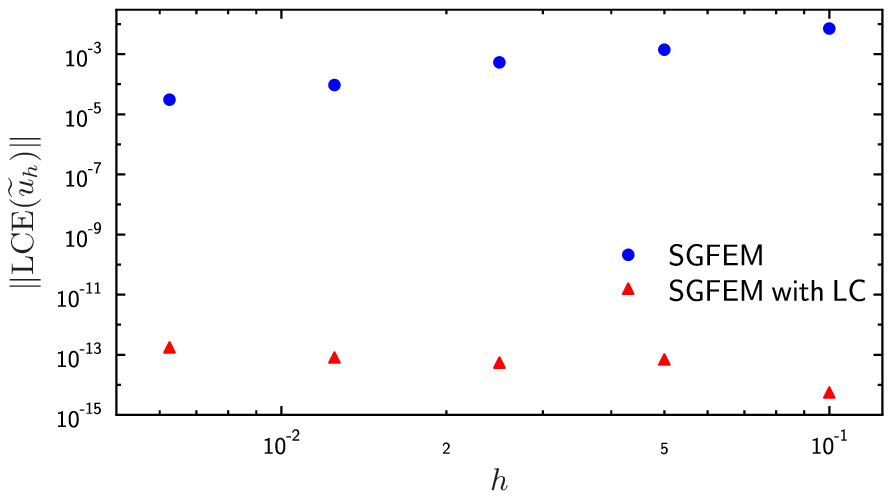}
\includegraphics[scale=0.65]{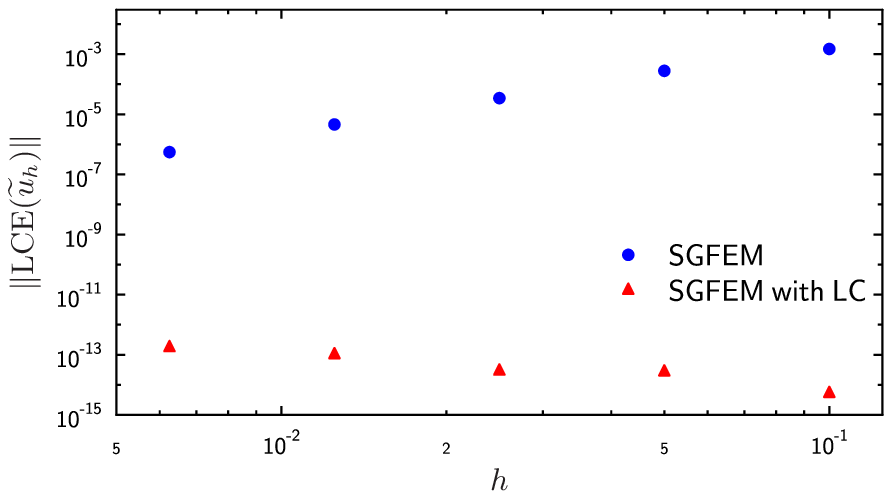}
\includegraphics[scale=0.65]{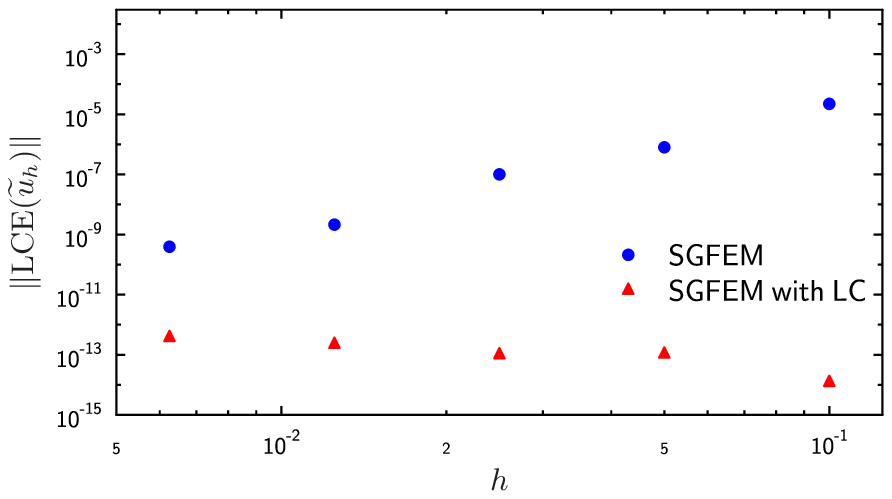}
\includegraphics[scale=0.65]{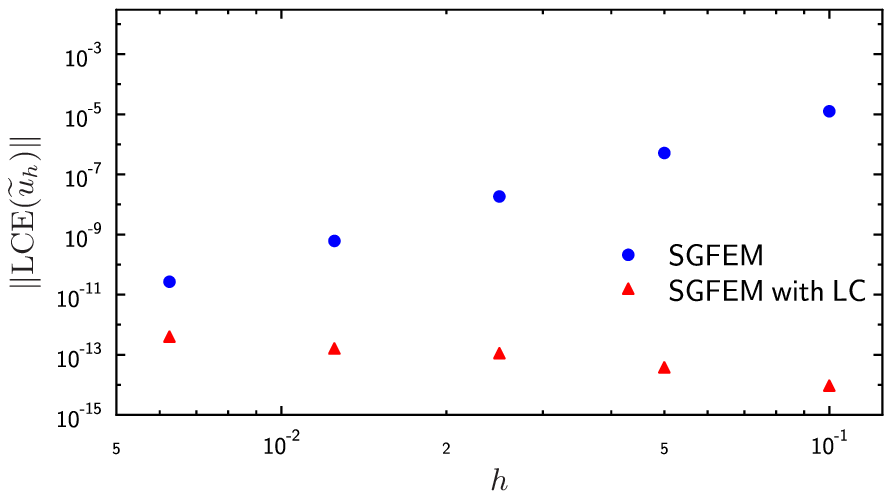}
\caption{Example 2: $\| \text{LCE}(\widetilde{u}_h) \|$ v.s. $h$ using  $p=1$ (top left), $p=2$ (top right), $p=3$ (bottom left) and $p=4$ (bottom right). SGFEM corresponds to $\widetilde{u}_h = u_{h,\tE} \in V_{h,\tE}^p$, and SGFEM with LC corresponds to $\widetilde{u}_h = \widehat{u}_{h,\tE} \in V_{h,\tE}^p$ obtained from \Cref{alg:3}. }
\label{figLCmeanerror}
\end{figure}

In \Cref{figLCmeanerror}, we also collect the mean absolute errors of the local conservation, which calculated as
$$
\| \text{LCE}(w) \|= \frac{1}{N^*}\sum_{\forall \tau^* \in \mathcal{T}^*} |\text{LCE}_{\tau^*} (w)|,
$$
for several values of $h$. Although $u_{h,\tE}$ does not satisfy the local conservation, it is evident that the mean absolute errors tend to decrease as $\mathcal{T}_h$ is refined. As for $\widetilde{u}_{h,\tE}$, the errors are less than $10^{-12}$ for any $h$, which is technically attributed to the errors of numerical integration and machine precision. Theoretically, these errors are equal to zero.

\begin{figure}[] 
\centering
\includegraphics[scale=0.65]{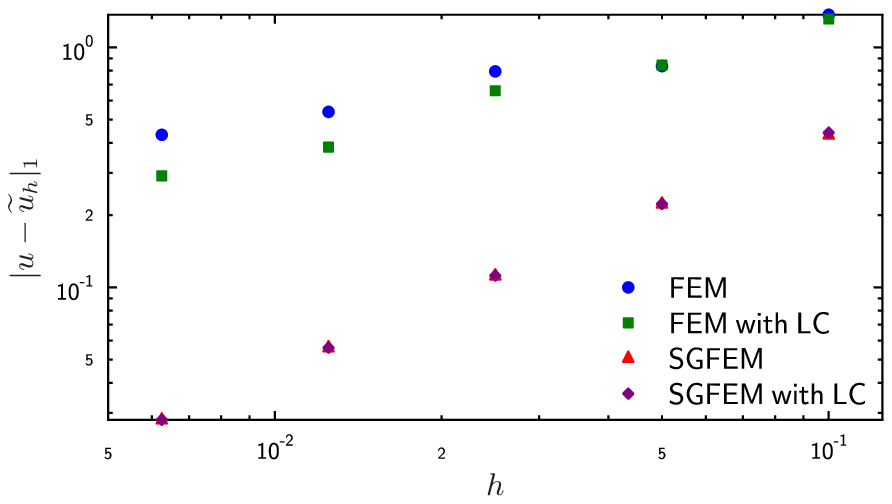}
\includegraphics[scale=0.65]{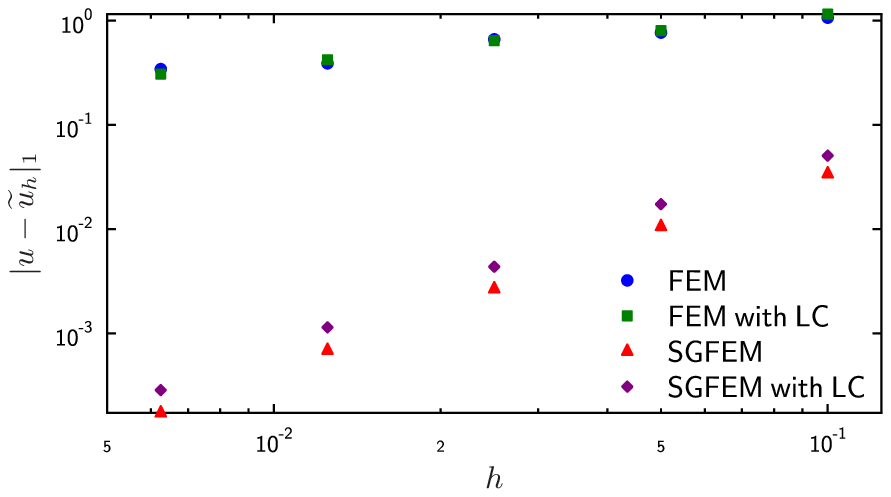}
\includegraphics[scale=0.65]{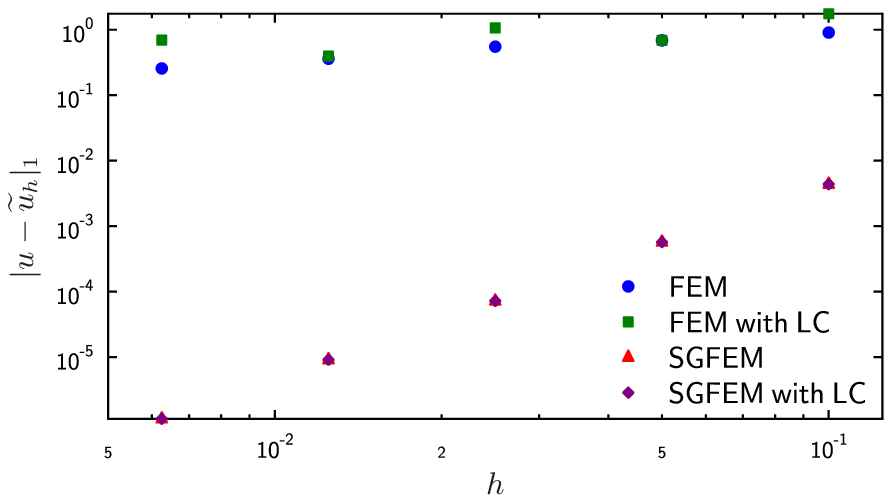}
\includegraphics[scale=0.65]{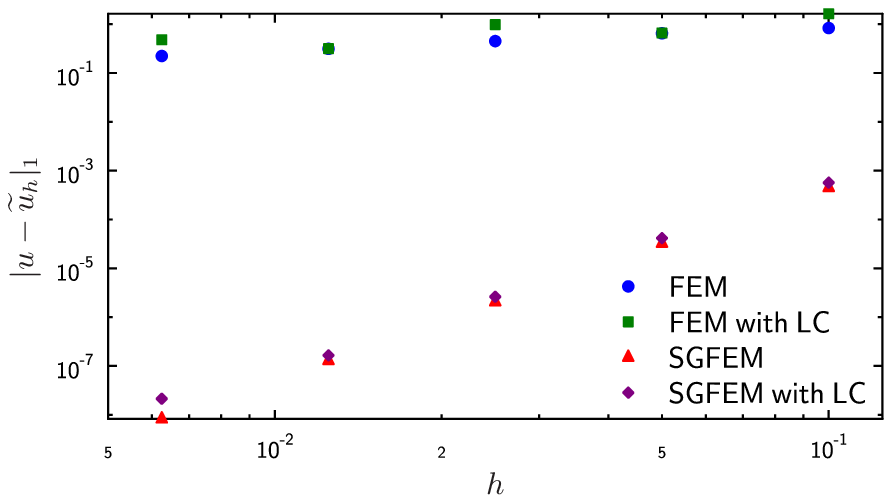}
\caption{Example 2: $|u - \widetilde{u}_h|_1$ v.s. $h$ using $p=1$ (top left), $p=2$ (top right), $p=3$ (bottom left), and $p=4$ (bottom right). FEM corresponds to $\widetilde{u}_h = u_h \in V_h^p$,
FEM with LC corresponds to $\widetilde{u}_h = \widehat{u}_h \in V_h^p$ obtained from \Cref{alg:3},
 SGFEM corresponds to $\widetilde{u}_h = u_{h,\tE} \in V_{h,\tE}^p$, and SGFEM with LC corresponds to $\widetilde{u}_h = \widehat{u}_{h,\tE} \in V_{h,\tE}^p$ obtained from \Cref{alg:3}.}
\label{figh1seminormLC}
\end{figure}

\begin{figure}[]
\centering
\includegraphics[scale = 0.65]{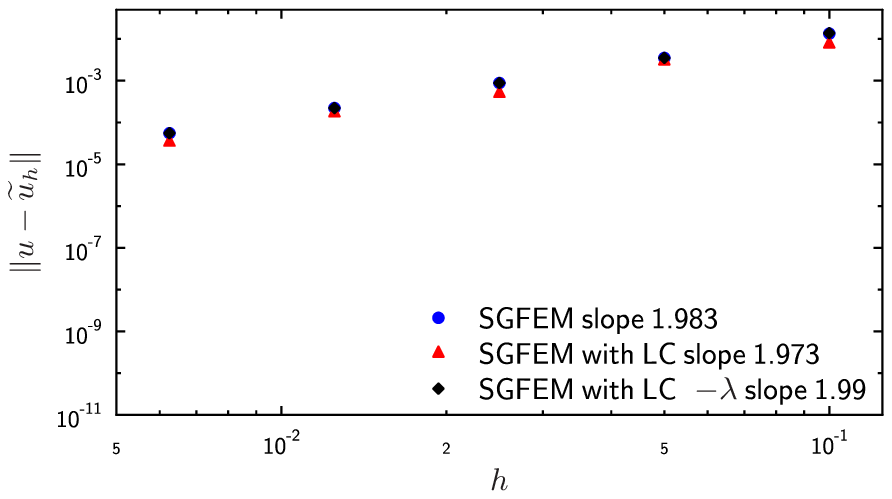}
\includegraphics[scale = 0.65]{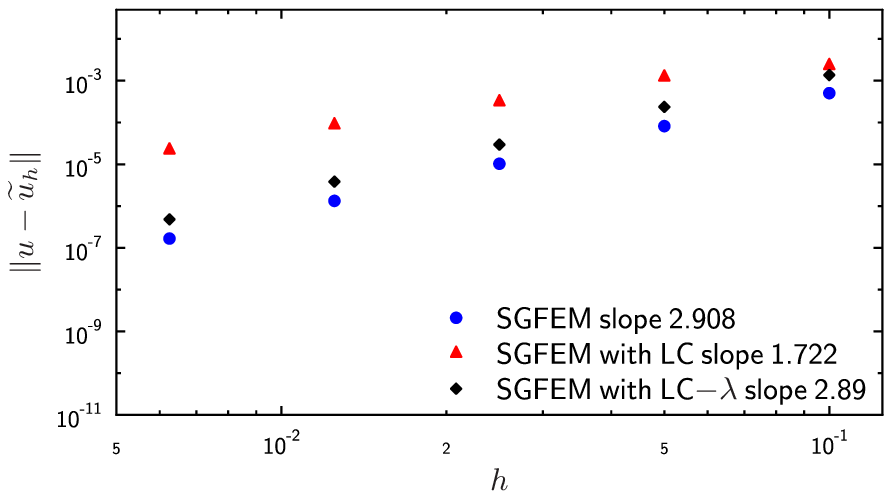}
\includegraphics[scale = 0.65]{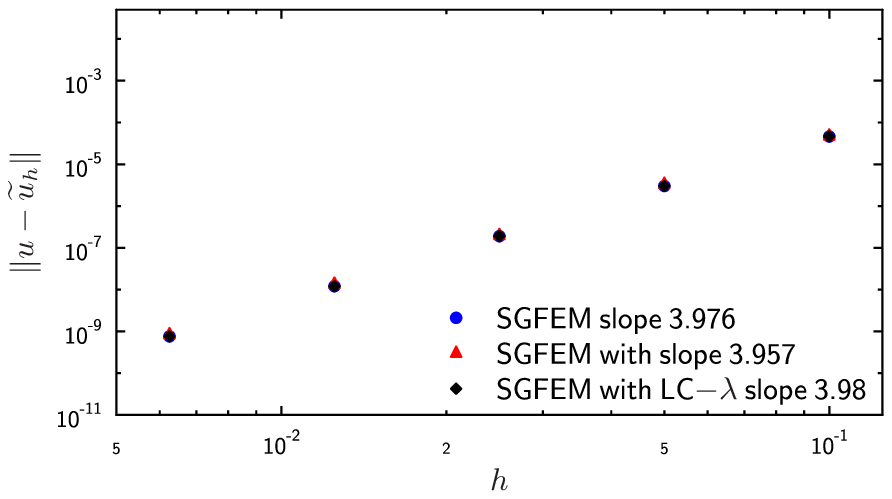}
\includegraphics[scale = 0.65]{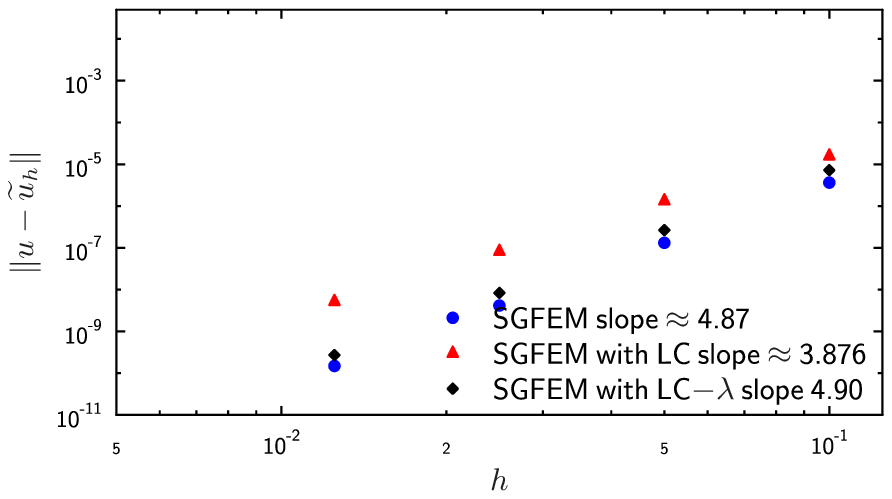}
\caption{Example 2: $\|u - \widetilde{u}_h\|$ v.s. $h$ using $p=1$ (top left) and $p=2$ (top right),
$p=3$ (bottom left), and $p=4$ (bottom right).
 SGFEM corresponds to $\widetilde{u}_h = u_{h,\tE} \in V_{h,\tE}^p$, and SGFEM with LC corresponds to $\widetilde{u}_h = \widehat{u}_{h,\tE} \in V_{h,\tE}^p$ obtained from \Cref{alg:3}, SGFEM with LC - $\lambda$ corresponds to $\widetilde{u}_h = \widehat{u}_{h,\tE} - \lambda$ obtained from \Cref{alg:3}.}
\label{figl2normLC}
\end{figure}

Comparison of $H^1$ semi-norm errors of FEM, SGFEM and the corresponding constrained problem for achieving local conservation can be seen in the \Cref{figh1seminormLC}. We can see that imposing the conservation in control volumes by Lagrange multipliers does not affect the optimal convergence. However, it is not the case for the corresponding errors in $L^2$ norm as shown in \Cref{figl2normLC}. This finding agrees with prior studies in  \cite{Abreu17}. As stated in this reference, the optimal convergence rate in $L^2$ norm can be recovered by adding the Lagrange multiplier as a corrector to the approximate solution, that is we calculate $\|u - \tilde{u}_h - \lambda\|_{L^2(\Omega)}$ where $\lambda$ is the Lagrange multiplier values over control volumes. 

\section{Conclusions} 
\setcounter{equation}{0}
\label{sec:CO}

In this paper, an application of SGFEM approximation for a  two-point boundary value problem whose elliptic coefficient is nonlinear and discontinuous has been presented. SGFEM relies on enriching the standard finite element space with auxiliary functions that capture the discontinuity effect. For problems with persistent discontinuity, SGFEM is preferred over standard continuous Galerkin finite element method. This is mainly due to its flexibility of not requiring the mesh to conform with the discontinuity configuration in order to maintain optimality of its convergence properties as affirmed by a detailed mathematical analysis of the method. A set of numerical examples with sufficiently high contrast coefficients verifies the aforementioned theoretical investigation.

As typical finite element approximations lack the local conservation property, a Lagrange multiplier technique is utilized to impose this property in each control volume as constraints to the minimization of an energy functional obtained from a linearization of the variational formulation.  An observation through numerical experiments reveals that using this technique, the optimal convergence property in $H^1$ semi-norm is preserved. An optimal convergence in $L^2$-norm can be recovered by using the Lagrange multiplier values as a corrector (see \cite{Abreu17}).

Extension of the present work includes an application of SGFEM approximation to Richards equation for unsaturated flow in heterogeneous soil. In this case, the soil heterogeneity is actualized by a layered system where every layer has its own constitutive relations, which is also dependent on the unknown function to be approximated. The matter is made more complicated due to the temporal dependence of Richards equation.
In addition to existing issues described in the present work, there is a need for an accurate and efficient time marching procedure that is suitably tailored to the generalized finite element methods.

\bibliography{mybibfile}   
\bibliographystyle{plain}  

\section*{Appendix A: Existence of the Solution to the Adjoint Problem and Its Approximation}
\setcounter{equation}{0}
\renewcommand{\thesection}{A} 
\counterwithin{proposition}{section} 
\numberwithin{equation}{section}

Let $u \in H^1_0(\Omega)$ be a weak solution of \eqref{eq:probstatment} and $\widetilde{u}_h \in \widetilde{V}_h \subset H^1_0(\Omega)$ be the corresponding approximate solution of $u$. 
Recall that $\widetilde{u}_h$ is partitioned into $N$ elements,  and the true solution $u$ is continuous and piecewise defined on $\Omega_j, j = 0,1,\cdots, m_\Gamma$. So  $\psi = u - \widetilde{u}_h$ is a continuous and piecewise function on $n = N+m_\Gamma$ subintervals in $\Omega$. Now let
$\mathcal{T}_j = \{ \tau \in  \mathcal{T}_h \setminus \mathcal{T}_{h,\Gamma} \text{ such that }  \tau \cap  \Omega_j \ne \varnothing  \}$ and  $\mathcal{T}_{j,\Gamma} = \{ \tau \cap \Omega_j \text{ such that } \tau \in \mathcal{T}_{h,\Gamma}  \}$
and set
$$\mathcal{S}_h = \{\varsigma_k = (s_{k-1},s_k), k=1,\cdots,n \text{ such that } \varsigma_k \text{ is either } \tau \in \mathcal{T}_j \text{ or } \tau_{1/2} \in \mathcal{T}_{j,\Gamma}, j= 1, \cdots, m_\Gamma \}.$$
Consider the following problem:
\begin{equation} \label{eq:probstatmentadj}
\begin{aligned}
&\text{find } \varphi \in H^1_0(\Omega) \text{ governed by}\\
&\begin{cases}
 \displaystyle - \frac{\d}{\d x}\left( \alpha_k(x) \frac{\d \varphi}{\d x} \right) + \beta_k(x) \frac{\d \varphi}{\d x} = \psi(x), ~~ x \in \varsigma_k, ~~ k = 1,\cdots, n,\\
 -\alpha_{k}(x)\dfrac{\d \varphi}{\d x} \Big|_{x = s_{k}^-} =  -\alpha_{k+1}(x)\dfrac{\d \varphi}{\d x} \Big|_{x = s_{k}^+} , ~~ k = 1,\cdots, n-1,\\
 \end{cases}
 \end{aligned}
 \end{equation}
where for $\sigma : [0,1] \to H^1_0(\Omega)$ defined by $\sigma(t) = \widetilde{u}_h + t \psi$,
$\alpha_k: \varsigma_k \cap \Omega_j \to \mathbb{R}$ and $\beta_k : \varsigma_k \cap \Omega_j \to \mathbb{R}$ are defined by  
\begin{equation*}
\begin{aligned}
\alpha_k(x) &= \int_0^1 \kappa_j(x, [\sigma(t)](x)) \, \d t, \\
\beta_k(x) &= \int_0^1 D_2 \kappa_j(x, [\sigma(t)](x)) [\sigma(t)]^\prime(x)\ \d t, \text{ with } [\sigma(t)]^\prime(x) = \widetilde{u}_h^\prime(x) + t \psi^\prime(x),\\
\end{aligned}
\end{equation*}
where $D_i\kappa_j$, $i =1,2$ implies the derivative of $\kappa_j$ with respect to the $i$th variable.  Since $u, \widetilde{u}_h \in H^1_0(\Omega)$ with $\Omega \subset \mathbb{R}$, $u, \widetilde{u}_h , \psi \in C(\overline{\Omega})$, which together with the assumption that $\kappa_j \in C^1(\overline{\Omega}_j \times \mathbb{R})$ implies that $\alpha_k \in C^1(\overline{\varsigma}_k)$. 
 Also, $\alpha_k(x) >0$, due to the positivity of $\kappa$, in particular we let $\kappa_{j,\min} \leq \alpha_k(x) \leq \kappa_{j,\max}$, and let $|D_i \kappa_j(x,v) | \leq \kappa'_{j,\max} <\infty$, for $(x,v) \in  \overline{\Omega}_j \times \mathbb{R}$. Also, set $\displaystyle \kappa_{\text{min}}= \min_{0\le j \le m_\Gamma} \kappa_{j,\text{min}}$ and $\displaystyle \kappa_{\text{max}} = \max_{0\le j \le m_\Gamma} \kappa_{j, \text{max}}$ and $\displaystyle \kappa'_{\text{max}} = \max_{0 \le j \le m_\Gamma} \kappa'_{j, \text{max}}$.
By this we can bound
\begin{equation} \label{eq:boundbeta}
|\beta_k(x)| \leq \kappa^\prime_{j,\max} \int_0^1 |u'(x) + t \psi'(x) | \ \d t \leq \kappa'_{\max} \big( |u'(x)| +  |\psi'(x)| \big), ~~ x \in \varsigma_k. 
\end{equation}

If there is $\varphi \in H^1_0(\Omega)$ satisfying \cref{eq:probstatmentadj}, then it also satisfies  \cref{eq:weakadjoint}. We establish the existence of a unique $\varphi$ by actually solving  \cref{eq:probstatmentadj} and exhibiting a formula for it. Due to the non-smoothness of the given data over $\Omega$, $\varphi$ is constructed  as a continuous function that is defined in a piecewise manner, which is denoted by $\varphi|_{\varsigma_k} = \varphi_k$. Using method of variation of parameters and integrating factor, it can be expressed as
\begin{align*}
\varphi_k(x) = a_{k}+ b_k I_k(x) + \int_{s_{k-1}}^x \dfrac{I_k(\xi)\psi(\xi)}{\eta_k(\xi)} \ \d \xi +  \int_x^{s_k} \dfrac{I_k(x) \psi(\xi) }{\eta_k(\xi)}\ \d \xi, ~~ x \in \varsigma_k, ~~k  = 1, \cdots, n,
\end{align*}
where $p_k: \varsigma_k \rightarrow \mathbb{R}$, $\eta_k: \varsigma_k \rightarrow \mathbb{R}$ and $I_k: \varsigma_k \rightarrow \mathbb{R}$ are defined as
\begin{equation*}
p_k(x) =  \int_{s_{k-1}}^x \dfrac{\beta_k(t)}{\alpha_k(t)}\ \d t,  \hspace{0.3cm}
\eta_k(x) = e^{p_k(x)}, \hspace{0.3cm} I_k(x) = \int_{s_{k-1}}^x \dfrac{\eta_k(t)}{\alpha_k(t)}\ \d t,
\end{equation*}
and $a_k, b_k$ are constants to be determined. The flux is
\begin{align*}
-\alpha_k(x)\varphi_k^\prime(x) = -\eta_k(x) \left( b_k + \int_x^{s_k} \dfrac{\psi(\xi)}{\eta_k(x)}\ \d \xi  \right),
\end{align*}
while the second derivative is
\begin{equation} \label{eq:phiprime}
\varphi_k^{\prime\prime}(x) = -\dfrac{\psi(x)}{\alpha_k(x)} + h_k(x) \left( b_k + \int_x^{s_k} \dfrac{\psi(\xi)}{\eta_k(\xi)}\ \d \xi \right),
\end{equation}
where
$$
h_k(x) =  \frac{\eta_k(x)(\beta_k(x) - \alpha_k'(x))}{(\alpha_k(x))^2}.
$$

Using \eqref{eq:boundbeta}, $p_k(x)$, $\eta_k(x)$, $I_k(x)$ for $x \in \overline{\varsigma}_k$ can be bounded as
\begin{equation*}
|p_k(x)| \leq  \frac{1}{\kappa_{\min}} \int_{s_{k-1}}^x  \hspace*{-0.3cm} |\beta_k(t)|\ \d t \leq 
\frac{\kappa'_{\max}}{\kappa_{\min}} \int_{\varsigma_k}  |u'(t) |+ |\psi'(t)| \d t\\
\le \rho_k,
\end{equation*}
\begin{equation} \label{eq:eta}
\begin{aligned}
 e^{-\rho_k}  \le \eta_k(x) &\le e^{\rho_k}, 
 \end{aligned}
\end{equation}

\begin{equation} \label{eq:Ik}
\begin{aligned}
\frac{|\varsigma_k| e^{-\rho_k}}{\kappa_{\max}} \le \frac{1}{\kappa_{\max}} \int_{s_{k-1}}^x \hspace*{-0.2cm} \eta_k(t) \d t \le
I_k(x) &\le \frac{1}{\kappa_{\min}} \int_{s_{k-1}}^x \hspace*{-0.2cm} \eta_k(t) \d t
\le \frac{|\varsigma_k|e^{\rho_k}}{\kappa_{\min}},
\end{aligned}
\end{equation}
where
\begin{equation*}
\rho_k = \frac{\sqrt{|\varsigma_k|}\kappa'_{\max}}{\kappa_{\min}} (|u|_{1,\varsigma_k} + |\psi|_{1,\varsigma_k}).
\end{equation*}
In some situations, we may globally bound $\rho_k$ as
\begin{equation*}
\rho_k \le \frac{\sqrt{|\Omega|}\kappa'_{\max}}{\kappa_{\min}} (|u|_1 + |\psi|_1) \le \frac{3\sqrt{|\Omega|}\kappa'_{\max} C_0^{-1} \| f \|}{\kappa_{\min}} =: \overline{\rho},
\end{equation*}
where we have used the boundedness of $u$ and $\widetilde{u}_h$ in terms of $f$ (see \Cref{thm:existence222} and its proof). Furthermore by Cauchy-Schwarz inequality,
\begin{equation*}
\begin{aligned}
\sum_{k=i}^j \rho_k &\le  \frac{\kappa'_{\max}}{\kappa_{\min}}  \Big( \sum_{k=i}^j |\varsigma_k|  \Big)^{1/2} \,
\left(  \Big(  \sum_{k=i}^j |u|^2_{1,\varsigma_k}  \Big)^{1/2}  +  \Big(  \sum_{k=i}^j |\psi|^2_{1,\varsigma_k}  \Big)^{1/2} \right)  \le  \overline{\rho}.
\end{aligned}
\end{equation*}

Straightforward calculation shows that, 
$$
\alpha_k'(x) = \displaystyle \int_0^1 D_1 \kappa_j(x,[\sigma(t)](x)) +  D_2 \kappa_j(x,[\sigma(t)](x)) [\sigma(t)]'(x) \ \d t,  
$$
from which 
$$
\begin{aligned}
|\beta_k(x) - \alpha'_k(x)|  &\le  \displaystyle \int_0^1 | D_1 \kappa_j(x,[\sigma(t)](x)) | \ \d t \leq  \kappa'_{\max},
\end{aligned}
$$
and thus
\begin{equation*}
\begin{aligned}
|h_k(x)| \leq \frac{\kappa'_{\max}e^{\rho_k}}{\kappa_{\min}^2} \le \frac{\kappa'_{\max}e^{\overline{\rho}}}{\kappa_{\min}^2}.
\end{aligned}
\end{equation*}

\setcounter{section}{0}
\setcounter{proposition}{0}
\renewcommand{\theproposition}{\thesection.\arabic{proposition}}
\renewcommand{\theequation}{\thesection.\arabic{equation}}

\begin{proposition} \label{pr:A1}
There exists a unique $\varphi \in H^1_0(\Omega)$ satisfying \cref{eq:probstatmentadj} such that
\begin{equation} \label{eq:h2varphi}
\sum_{k = 1}^n| \varphi |^2_{2,\varsigma_k}  \le C \| \psi \|^2,
\end{equation}
where $C$ depends only on $\kappa_{\min}$, $\kappa_{\max}$, $\kappa'_{\max}$, $|\Omega|$, and $\|f\|$.
\end{proposition}
\begin{proof}
First of all, we need to find the collection of $\{a_k\}$ and $\{b_k\}$ which will be determined by imposing the following conditions:
\begin{itemize}
\item[-]  The boundary conditions, we get 2 equations:
 \begin{align*}
 a_1 &= 0
 \\
 a_{n} + b_{n} I_{n}(s_{n}) &= - \int_{\varsigma_{n}} \frac{I_{n}(\xi) \psi(\xi)}{\eta_{n}(\xi)} \ \d \xi.
 \end{align*}
\item[-]  The continuity of the function solution:  $\varphi_i(s_i) = \varphi_{i+1}(s_i)$,  we get $n-1$ equations:
 \begin{equation*}
a_i +  b_i I_i(s_i) - a_{i+1} = - \int_{\varsigma_i} \frac{I_i(\xi)\psi(\xi)}{\eta_i(\xi)} \ \d \xi,
\hspace{0.2cm} i=1, \cdots, n - 1.
\end{equation*}
\item[-] The continuity of the flux: $-\alpha_i(s_i) \varphi_i^\prime(s_i) = -\alpha_{i+1}(s_i) \varphi_{i+1}^\prime(s_i)$, we also get $n - 1$ equations:
\begin{equation*}
b_{i} \eta_{i}(s_i) - b_{i+1}  =  \int_{\varsigma_{i+1}} \frac{\psi(\xi)}{\eta_{i+1}(\xi)} \ \d \xi, 
\hspace{0.2cm} i=1, \cdots, n - 1.
\end{equation*}
\end{itemize} 

Altogether we obtain a square matrix system sized $2n$ that can be written as,
{\small
\begin{equation} \label{eq:bmatrix}
\left[ \begin{array}{@{}ccccc|ccccc@{}}
1 & 0  & 0 & \ldots & 0 & 0  & 0 & \ldots & 0 & 0\\
1 & -1 & 0 & \ldots  & 0 & I_1(s_1) & 0 & \ldots & 0 & 0 \\
0 & 1 & -1 & \ldots  & 0 & 0 & I_2(s_2) & \ldots & 0 & 0 \\
\vdots &  & \ddots & \ddots  & \vdots& \vdots &  & \ddots & & \vdots \\
0 & 0 & & 1 & -1 & 0 & 0 & & I_{n-1}(s_{n-1}) & 0 \\ \hline
0 & 0 &\ldots & 0 & 0 & \eta_1(s_1) & -1 & 0 & \ldots & 0 \\
0 & 0 & \ldots & 0 & 0 & 0 & \eta_2(s_2) & -1 & \ldots & 0\\
\vdots&  & \ddots & & \vdots & \vdots & & \ddots & \ddots & \vdots \\
0 & & \ldots & & 0 &0 & \ldots & & \eta_{n-1}(s_{n-1}) & -1\\
0 & & \ldots & 0 & 1 & 0 & \ldots & & 0 & I_n(s_n)
\end{array} \right]
\left[\begin{array}{c}
a_1\\
a_2\\
a_3\\
\vdots\\
a_n\\ \hline
b_1\\
b_2 \\
\vdots\\
b_{n-1}\\
b_n
\end{array}\right]
=
\left[\begin{array}{c}
F_1\\
F_2\\
F_3\\
\vdots\\
F_n\\ \hline
G_1\\
G_2 \\
\vdots\\
G_{n-1}\\
G_n
\end{array}\right],
\end{equation}}
where 
\begin{equation*}
\begin{aligned}
F_ 1 &= 0, ~~F_{i+1}= -\int_{\varsigma_{i}} \hspace{-0.1cm}\frac{I_{i}(\xi) \psi(\xi)}{\eta_{i}(\xi)}\d \xi , ~ G_i  = \int_{\varsigma_{i+1}}\hspace{-0.1cm} \frac{\psi(\xi)}{\eta_{i+1}(\xi)} \d \xi ,~~ G_{n} = -\int_{\varsigma_{n}}\hspace{-0.2cm} \frac{I_n(\xi) \psi(\xi)}{\eta_{n}(\xi)} \d \xi,
\end{aligned}
\end{equation*}
for $ i = 1, \cdots, n-1$.
By \eqref{eq:eta}, \eqref{eq:Ik} and Cauchy-Schwarz inequality,
\begin{equation*}
\begin{aligned}
|F_{i+1}| &\le \int_{\varsigma_{i}} \Big| \dfrac{I_{i}(\xi)\psi(\xi)}{\eta_{i}(\xi)} \Big | \d \xi 
\le \frac{|\varsigma_{i}| e^{2 \rho_{i}}
 }{\kappa_{\min}} \int_{\varsigma_{i}} \hspace{-0.05cm} |\psi(\xi)| \d \xi
\le \frac{ e^{2\overline{\rho}}\ |\varsigma_{i}|^{3/2} \|\psi \|_{0,\varsigma_{i}}}{\kappa_{\min}} .
\end{aligned}
\end{equation*}
Also $G_i$ can be bounded using \eqref{eq:eta} and Cauchy - Schwarz inequality,
\begin{equation} \label{eq:Gi}
\begin{aligned}
|G_i| &\le \int_{\varsigma_{i+1}} \Big| \dfrac{\psi(\xi)}{\eta_{i+1}(\xi)} \Big| \d \xi \le e^{\rho_{i+1}}
\int_{\varsigma_{i+1}} \hspace{-0.3cm}|\psi(\xi)| \ \d \xi
\le e^{\overline{\rho}} \ |\varsigma_{i+1}|^{1/2} \|\psi\|_{0,\varsigma_{i+1}} .
\end{aligned}
\end{equation}
In a similar fashion to estimating $F_{i+1}$,
\begin{equation} \label{eq:boundgn}
\begin{aligned}
|G_n| \le  \frac{e^{2\overline{\rho}}\ |\varsigma_{n}|^{3/2} \|\psi \|_{0,\varsigma_{n}}}{\kappa_{\min}} .\end{aligned}
\end{equation}

Together with \cref{eq:eta} and \cref{eq:Ik}, all coefficients in the system \cref{eq:bmatrix} are shown to be bounded.  By the standard row reduction, we can transform the matrix in \cref{eq:bmatrix} into an upper triangular matrix with nonzero diagonal entries. Using back substitution, there is a unique collection of  constants $\{a_k\}$ and $\{b_k\}$ such that $\varphi$ is the unique solution of the adjoint problem.
In particular,
$$ 
b_i = b_n \displaystyle \prod_{j = i}^{n-1}(\eta_j(s_j))^{-1} + \sum_{l = i}^{n-1} G_l  \displaystyle \prod_{j = i}^{l}( \eta_j(s_j))^{-1},~~ i = 1, \cdots, n-1,
$$
and %
$$
b_n = \dfrac{\eta_{n-1}(s_{n-1})(G_n + F_n ) - G_{n-1}I_{n-1}(s_{n-1})}{I_n(s_n)\eta_{n-1}(s_{n-1}) + I_{n-1}(s_{n-1})}.
$$
Notice that $\varphi \in H^1_0(\Omega)$, which is achieved from imposing continuity of $\varphi$ and its flux on every $s_k$.

At this stage we perform various estimations to bound $\{ b_i \}$.
By taking the absolute value on $b_n$, we can estimate
\begin{equation*}
\begin{aligned}
|b_n| &\leq \frac{ |G_n| }{I_n(s_n)} + \dfrac{\eta_{n-1}(s_{n-1})|F_n|}{I_{n-1}(s_{n-1})} + |G_{n-1}|.
\end{aligned}
\end{equation*}
Using the bound of $G_n$ in \eqref{eq:boundgn} and the lower bound of $I_n$ in \cref{eq:Ik},  the first term is bounded by 
\begin{equation*}
\begin{aligned}
\frac{ |G_n| }{I_n(s_n)} \leq \frac{\kappa_{\max} e^{3\overline{\rho}}\ |\varsigma_{n}|^{1/2}  \|\psi\|_{0,\varsigma_{n}}}{\kappa_{\min}},
\end{aligned}
\end{equation*}
while the second term is bounded by 
\begin{equation*}
\begin{aligned}
 \dfrac{\eta_{n-1}(s_{n-1})|F_n|}{I_{n-1}(s_{n-1})} \leq \frac{ \kappa_{\max} e^{4 \overline{\rho}} |\varsigma_{n-1}|^{1/2}  \|\psi \|_{0,\varsigma_{n-1}}}{\kappa_{\min}}.
\end{aligned}
\end{equation*}
Combining these last two estimates and using \cref{eq:Gi} gives
\begin{equation} \label{eq:bn}
|b_n| \le \dfrac{\kappa_{\max} e^{4 \overline{\rho}} }{ \kappa_{\min}} (2 |\varsigma_{n}|^{1/2}  \|\psi\|_{0,\varsigma_{n}}  + |\varsigma_{n-1}|^{1/2}  \|\psi \|_{0,\varsigma_{n-1}}  ).
\end{equation}
Using the bound of $G_i$ and $b_n$ and Cauchy-Schwarz inequality, 
\begin{equation*}
\begin{aligned}
|b_n| + \sum_{l =i}^{n-1}  |G_l|  
\le  \dfrac{3 \kappa_{\max} e^{4 \overline{\rho}} }{ \kappa_{\min}} \sum_{l= 1}^{n}|\varsigma_{l}|^{1/2}  \|\psi\|_{0,\varsigma_{l}} 
\le \dfrac{3 \kappa_{\max} e^{4 \overline{\rho}} }{ \kappa_{\min}}  |\Omega|^{1/2} \|\psi\|.
\end{aligned}
\end{equation*}
Now for $i = 1, \cdots, n-1$, we have
\begin{equation*}
\begin{aligned}
|b_i| & \le |b_n| \exp\left(\sum_{j=i}^{n-1} \rho_j\right) + \sum_{l =i}^{n-1}  |G_l | \exp\left( \sum_{j = i}^l \rho_j\right)
\\
&\le \exp\left( \sum_{j=i}^{n-1} \rho_j\right) \Big(|b_n| +  \sum_{l =i}^{n-1}  |G_l|\Big)
\\
&\le  \dfrac{3 \kappa_{\max} e^{5 \overline{\rho}} }{ \kappa_{\min}} |\Omega|^{1/2}  \|\psi\|.
\end{aligned}
\end{equation*}
Thus we can write that 
$$
|b_k| \leq \dfrac{3 \kappa_{\max} e^{5 \overline{\rho}} }{ \kappa_{\min}}|\Omega|^{1/2} \|\psi\|, ~~ k = 1, \cdots ,n.
$$

Next, we take the squared power of \eqref{eq:phiprime} to get
\begin{equation*}
\begin{aligned}
|\varphi''_k(x)|^2 = \left(-\dfrac{\psi(x)}{\alpha_k(x)} + b_k h_k(x) + h_k(x) \int_x^{s_k} \frac{\psi(\xi)}{\eta_k(\xi)} \d \xi \right)^2 \leq 3(J_1(x) + J_2(x) + J_3(x)),
\end{aligned}
\end{equation*}
where
\begin{equation*}
\begin{aligned}
J_1(x) &= \Big| \dfrac{\psi(x)}{\alpha_k(x)}\Big|^2, \hspace{0.3cm} J_2(x) = |b_k h_k(x)|^2, \hspace{0.3cm}  J_3(x) = \left| h_k(x) \int_x^{s_k}\dfrac{\psi(\xi)}{\eta_k(\xi)}\d \xi\right|^2.
\end{aligned}
\end{equation*}
By integrating $J_1$, $J_2$, and $J_3$ over $\varsigma_k$ and sum the results up over $k$ from 1 to $n$, we get 
\begin{equation*}
\begin{aligned}
\sum_{k=1}^n \int_{\varsigma_k} J_1(x)\d x  & \le  \dfrac{1}{\kappa_{\min}^2} \int_\Omega |\psi(x)|^2 \ \d x  =  \dfrac{1}{\kappa_{\min}^2}  \|\psi\|^2,
\\ \\
\sum_{k=1}^n \int_{\varsigma_k} J_2(x) \d x & \le \left(\dfrac{3 e^{6 \overline{\rho}}  \kappa_{\max} \kappa'_{\max}}{ \kappa_{\min}^2}|\Omega|^{1/2} \|\psi\|  \right)^2 \int_{\Omega} 1 \d x 
= \left(\dfrac{3  e^{6 \overline{\rho}}  \kappa_{\max}\kappa'_{\max}}{ \kappa_{\min}^2} \right)^2 |\Omega|^2 \|\psi\|^2,
\\ \\
\sum_{k=1}^n \int_{\varsigma_k} J_3(x) \d x &\le   
\left(\frac{\kappa'_{\max} e^{\overline{\rho}}}{\kappa_{\min}}\sum_{k = 1}^n \int_{\varsigma_k} \psi_k(\xi) \d x \right)^2 
\le \left(\frac{\kappa'_{\max}e^{\overline{\rho}} }{\kappa_{\min} } \right)^2 |\Omega| \|\psi\|^2.
\end{aligned}
\end{equation*}
Thus
\begin{equation*}
\begin{aligned}
\sum_{k=1}^n |\varphi|_{2,\varsigma_k}^2 \le C \|\psi\|^2,
\end{aligned}
\end{equation*}
where
$$
C = \dfrac{1}{\kappa_{\min}^2} + \left(\dfrac{3  e^{6 \overline{\rho}}  \kappa_{\max}\kappa'_{\max}}{ \kappa_{\min}^2} \right)^2 |\Omega|^2 +  \left(\frac{\kappa'_{\max}e^{\overline{\rho}} }{\kappa_{\min} } \right)^2 |\Omega|. 
$$

\end{proof}

\begin{remark}
If $\kappa_j(x,u(x)) = \kappa_j(u(x))$, then
\begin{equation*}
\begin{aligned}
\alpha_k(x) &= \int_0^1 \kappa_j([\sigma(t)](x)) \, \d t, \\
\beta_k(x) &= \int_0^1  \kappa_j'( [\sigma(t)](x)) [\sigma(t)]^\prime(x)\ \d t, \text{ with } [\sigma(t)]^\prime(x) = \widetilde{u}_h^\prime(x) + t \psi^\prime(x),\\
\end{aligned}
\end{equation*}
which yields
$\beta_k(x) - \alpha'_k(x) =0,$ and thus $h_k(x) = 0$, and  \eqref{eq:phiprime} becomes
$$
\varphi_k^{\prime\prime}(x) = -\dfrac{\psi_k(x)}{\alpha_k(x)}. 
$$
By squaring both sides and integrate over $\varsigma_k$, we can estimate
\begin{equation*}
\begin{aligned}
|\varphi_k|^2_{2,\varsigma_k} = \int_{\varsigma_k} |\varphi''_k|^2 \d x = \int_{\varsigma_k} \left| \frac{\psi_k(x)}{\alpha_k(x)}\right|^2 \d x \leq (\kappa_{\min})^{-2} \|\psi\|^2_{0,\varsigma_k}.
\end{aligned}
\end{equation*}
By summing up over $k$ we obtain
$$
\sum_{k=1}^n |\varphi|^2_{2,\varsigma_k} \leq (\kappa_{\min})^{-2} \|\psi\|^2.
$$
\end{remark}

\begin{proposition} \label{pr:A2}
Let $\varphi$ be the solution of the adjoint problem in \cref{eq:probstatmentadj},  and $ \mathcal{I}_{h,{\emph \tE}}^1 \varphi \in V_{h,{\emph \tE}}^1$ be its interpolant. Then,
$$
| \varphi - \mathcal{I}_{h,\tE}^1 \varphi |_1 \leq C h \| \psi \|.
$$
\end{proposition}
\begin{proof}
The proof is rather similar to the proof of \Cref{lem:norm16}:
\begin{align*}
|\varphi - \mathcal{I}_{h,\tE}^1 \varphi|^2_{1,\Omega_j} &= \sum_{\tau \in \mathcal{T}_j} | \varphi - \mathcal{I}_{\tau}^1 \varphi|^2_{1,\tau} + \sum_{\tau_{1/2} \in \mathcal{T}_{j,\Gamma}}  | \varphi - \mathcal{I}_{\tau,\tE}^1 \varphi|^2_{1,\tau_{1/2}}\\
&\le \sum_{\tau \in \mathcal{T}_j} h^2 \|\varphi''\|_{0,\tau}^2 + \sum_{\tau_{1/2} \in \mathcal{T}_{j,\Gamma}} h^2 \|\varphi''\|_{0,\tau_{1/2}}^2. 
\end{align*}
Summing up over every $\Omega_j$ and using \Cref{pr:A1} give the desired estimate.
\end{proof}

\section*{Appendix B: Some Analytical Solutions}
\setcounter{equation}{0}
\renewcommand{\thesection}{B} 

\begin{enumerate}
\item With elliptic coefficient \cref{eq:prob1}, the analytical solution of example 1 is 
\begin{equation} \label{eq:anasol1}
u(x) = 
\begin{cases}
\dfrac{1}{a_0} \ln(a_0(-5x^3/6 + C_1x) + 1), &x\in(0, \gamma_1],\
 \vspace{0.2cm}
 \\
\dfrac{1}{a_1} \ln(a_1(-5x^3/6 + C_1x + C_2)), & x \in (\gamma_1, \gamma_2], 
 \vspace{0.2cm}
 \\
\dfrac{1}{a_2} \ln(a_2 5(1-x^3)/6 + a_2C_1(x -1) x + 1), & x \in(\gamma_2,1),
\end{cases}
\end{equation} 
where $C_1$ and $C_2$ are solutions of the following system of nonlinear equations:  
\begin{equation*}
\begin{cases}
 \dfrac{1}{a_0} \ln (a_0(-5\gamma_1^3/6 +C_1\gamma_1 )+ 1) - 
 \dfrac{1}{a_1} \ln (a_1(-5\gamma_1^3/6 +C_1\gamma_1 + C_2)) = 0
 \vspace{0.2cm}
\\
\dfrac{1}{a_1} \ln (a_1(-5\gamma_2^3/6 + C_1\gamma_2 + C_2)) -
 \dfrac{1}{a_2} \ln (a_2 5(1 -\gamma_2^3)/6  + a_2 C_1(\gamma_2 - 1) + 1) =0,
 \end{cases}
\end{equation*}
which can be solved by a standard nonlinear solver.
\\

\item With elliptic coefficient \cref{eq:prob2}, the analytical solution of example 2 is 

\begin{equation} \label{eq:anasol2}
u(x) = 
\begin{cases}
-\ln\Big(- \dfrac{\sin(\pi x)}{a_0 \pi^2}  +\dfrac{x}{a_0 \pi} - C_1x +1 \Big), & x \in (0,\gamma_1],
\vspace{0.3cm}
\\ 
-\ln\Big(\dfrac{\sin(\pi\gamma_1) - \sin(\pi x)}{a_1\pi^2} + (1 + a_3 \pi C_2)\dfrac{( \gamma_1- x) }{a_1\pi} + C_3, & x\in (\gamma_1,\gamma_2],
\vspace{0.3cm}
\\
-\ln\Big(\dfrac{\sin(\pi\gamma_2) - \sin(\pi x)}{a_2\pi^2} + (1 + a_3\pi C_2)\dfrac{( \gamma_2 - x)}{a_2\pi} + C_4, & x \in(\gamma_2,\gamma_3],
\vspace{0.3cm}
\\
-\ln\Big(\dfrac{-\sin(\pi x)}{a_3 \pi ^2} + \dfrac{1-x}{a_3 \pi} + C_2(1-x)+ 1 \Big), &x\in(\gamma_3,1),
\end{cases}
\end{equation}
where 
\begin{equation*}
\begin{aligned}
C_2 &= \dfrac{p + q}{r}, ~~
p =  \dfrac{\sin(\pi \gamma_3)}{a_3 \pi^2} + \dfrac{\sin(\pi \gamma_2) - \sin(\pi\gamma_3) }{a_2 \pi^2}+ \dfrac{\sin(\pi \gamma_1) - \sin(\pi \gamma_2)}{a_1\pi^2} -\dfrac{ \sin(\pi \gamma_1)}{a_0\pi^2},
\\
q &= -\dfrac{\gamma_1}{a_0\pi} + \dfrac{\gamma_1 - \gamma_2}{a_1 \pi} + \dfrac{\gamma_2 - \gamma_3}{a_2\pi} +\dfrac{\gamma_3 - 1}{a_3\pi},\\
r &= 1 - \gamma_3 + \dfrac{a_3\gamma_1}{a_0} +\dfrac{a_3(\gamma_2 - \gamma_1)}{a_1} + \dfrac{a_3(\gamma_3 - \gamma_2)}{a_2} ,
\\
C_1 &= \dfrac{2}{a_0 \pi} + \dfrac{a_3 C_2}{a_0}, ~~
C_3 = \dfrac{-\sin(\pi\gamma_1)}{a_0 \pi^2} + \dfrac{\gamma_1}{a_0 \pi} + 1 - \gamma_1 C_1,
\\
C_4  &= \dfrac{\sin(\pi \gamma_1) - \sin(\pi \gamma_2)}{a_1 \pi^2} + \Big(  \dfrac{1}{a_1\pi} + \dfrac{a_3C_2}{a_1}\Big)(\gamma_1 - \gamma_2) + C_3.
\end{aligned}
\end{equation*}
\end{enumerate}


%
%


\end{document}